\documentclass[10pt]{article}
\usepackage[utf8]{inputenc}
\usepackage[letterpaper,top=1in,bottom=1in,left=1.2in,right=1.2in]{geometry}
\usepackage[english]{babel}
\usepackage{amsmath}
\usepackage{amsthm}
\usepackage{bm}
\usepackage{amssymb}
\usepackage{dsfont}
\usepackage{commath}
\usepackage{mathtools}
\usepackage{float}
\usepackage{algorithm}
\usepackage{algpseudocode}

\usepackage[numbers]{natbib}

\usepackage{xcolor}

\def\eps{\epsilon}
\def\be{\begin{equation}}
\def\ee{\end{equation}}
\def\bea{\begin{align}}
\def\eea{\end{align}}
\def\bea*{\begin{align*}}
\def\eea*{\end{align*}}
\def\bsx{\boldsymbol{\xi}}

\newtheorem{theorem}{Theorem}[section]
\newtheorem*{theorem*}{Theorem}

\newtheorem{lemma}{Lemma}[section]
\newtheorem{corollary}{Corollary}[section]
\newtheorem{proposition}{Proposition}[section]
\newtheorem{assumption}{Assumption}[section]

\theoremstyle{definition} 

\newtheorem{example}{Example}
\newtheorem{remark}{Remark}[section]

\DeclareMathOperator{\E}{\mathbb{E}}
\DeclareMathOperator{\F}{\mathbb{F}}

\DeclareMathOperator{\N}{\mathbb{N}}

\DeclareMathOperator{\R}{\mathbb{R}}

\DeclareMathOperator{\T}{\mathbb{T}}

\DeclareMathOperator{\Z}{\mathbb{Z}}

\DeclareMathOperator{\cA}{\mathcal{A}}
\DeclareMathOperator{\cB}{\mathcal{B}}
\DeclareMathOperator{\cC}{\mathcal{C}}

\DeclareMathOperator{\cE}{\mathcal{E}}
\DeclareMathOperator{\cF}{\mathcal{F}}
\DeclareMathOperator{\cG}{\mathcal{G}}

\DeclareMathOperator{\cI}{\mathcal{I}}
\DeclareMathOperator{\cJ}{\mathcal{J}}

\DeclareMathOperator{\cL}{\mathcal{L}}
\DeclareMathOperator{\cM}{\mathcal{M}}
\DeclareMathOperator{\cN}{\mathcal{N}}

\DeclareMathOperator{\cP}{\mathcal{P}}

\DeclareMathOperator{\cT}{\mathcal{T}}

\DeclareMathOperator{\cW}{\mathcal{W}}
\DeclareMathOperator{\cX}{\mathcal{X}}
\DeclareMathOperator{\cY}{\mathcal{Y}}
\DeclareMathOperator{\cZ}{\mathcal{Z}}

\DeclareMathOperator{\mathfrakC}{\mathfrak{C}}

\DeclareMathOperator{\mathfrakJ}{\mathfrak{J}}

\usepackage{subcaption}
\usepackage{graphicx}
\usepackage{color}

\usepackage{xcolor}
\usepackage{hyperref}
\hypersetup{colorlinks,linkcolor=blue,citecolor=blue,urlcolor = violet}
\usepackage[capitalise]{cleveref}

\numberwithin{equation}{section}
\usepackage{fancyhdr}
\title{Learning algorithms for mean field optimal control
\footnote{Research of HMS and QY was
partially supported by the National Science Foundation grant DMS 2406762.}}

\author{H. Mete Soner\footnote{Department of Operations Research and Financial
Engineering, Princeton University, Princeton, NJ, 08540, USA, email: 
{\tt soner@princeton.edu}. }
\and Josef Teichmann \footnote{Department of Mathematics, ETH, Z\"urich, Switzerland, 
email: {\tt jteichma@math.ethz.ch}.}
\and Qinxin Yan\footnote{Program in Applied and Computational
Mathematics, Princeton University, Princeton, NJ, 08540, USA, email: 
{\tt qy3953@princeton.edu}. }}

\date{\today}
\usepackage{oubraces}
\begin{document}
\maketitle

\vspace{0mm}
\begin{abstract}
We analyze an algorithm to numerically solve
mean-field optimal control problems by approximating the optimal feedback controls 
using neural networks with problem specific architectures.
We approximate the model
by an $N$-particle system and leverage
the exchangeability of the particles
to obtain substantial computational
efficiency.
In addition to several
numerical examples,
a convergence analysis 
is provided. We also develop a universal approximation theorem on Wasserstein spaces.
\end{abstract}
\vspace{1mm}

\vspace{3pt}

\def\d{\mathrm{d}}

\section{Introduction}
\label{sec:intro}

Mean field control (MFC) and mean field games (MFG) have emerged over the past two decades as important frameworks for modeling, analyzing and understanding large-population systems.  The subject originated from pioneering works of Lasry and Lions \cite{lasry_jeux_2006,lasry_jeux_2006-1,lasry_mean_2007} and independently of Huang, Caines, and Malham\'e \cite{huang_individual_2003,huang_large_2006,huang_nash_2007}, who introduced a new paradigm that overcomes the curse of dimensionality in classical $N$-particle differential games.  Subsequently, the
classical book of Carmona and Delarue \cite{CD} develops the
probabilistic approach and provides an 
excellent summary of the impressive theoretical foundations and the scope of MFG.  
These mathematical formulations are receiving increasing attention in various application domains, including finance, engineering, and social sciences.

Despite their conceptual elegance, MFG and MFC often lead to partial differential equations in distribution spaces or forward-backward stochastic differential equations (FBSDEs) with non-linear distributional dependences.  Numerically solving such equations is notoriously challenging, and many recent developments have focused on machine learning (ML) and deep neural networks (DNN) to alleviate these difficulties.  Classical PDE-based methods (see, e.g., \cite{reisinger,osher}) are suitable when the dimensionality remains moderate or when specific structural assumptions allow for tractable discretization schemes. However, in truly high-dimensional or more flexible settings, ML-based methods have emerged as practical and powerful alternatives.  

\medskip
\noindent\emph{Existing ML-based approaches for mean field control.}
Machine learning algorithms for solving high-dimensional control problems, including mean field control, have undergone intensive developments in recent years. One of the most direct approaches is to first approximate the mean field problem by a system of $N$ interacting particles, then solve the resulting finite-dimensional stochastic control problem by employing DNN-based methods such as those introduced in \cite{PhamConv,PhamNum,becker2019deep,DH,HanE,HanJentzenE,RS,RST}
and for partial differential equations in \cite{beck2020overview,grohs2023proof}.  

In particular, Pham and Warin \cite{PW2,PW1} employ deep neural networks to compute the optimal control of the $N$-particle approximation by simulating \emph{many independent $N$-particle systems} to ensure approximation accuracy.  By contrast, our proposed algorithm leverages the exchangeability of the particles and only requires \emph{one simulation} of the large-population system, which yields substantial computational speed-ups while still preserving approximation quality.  In that sense, our approach differs significantly from \cite{PW2,PW1} and we provide a comparison in the example
of systematic risk model in subsection \ref{ssec:risk}.   Furthermore, reinforcement-learning methods for learning the model along with the optimal policy have been explored in \cite{PhamActorCritic}, but we do not consider such model-learning features.

In another closely related  series of papers, Carmona and Lauri\`ere \cite{CL1,CL2,CL3} use deep learning to solve MFG or mean field control problems under a large-population viewpoint.  
They take a single large $N$ simulation for each epoch and train feedback policies 
that depend only on time and the individual state of each particle, but not on the distribution.  
In many scenarios, this yields effective approximations to the optimal MFC solution.  

The distributional dependence is included 
by Dayan\i kl\i, Lauri\`ere and Zhang~\cite{DLZ}  who use
the symmetry of the problem to construct efficient network architectures akin to \cite{GLPW}.  
Our approach also exploits this symmetry. However, instead of treating the $N$-particle system as a single large system and training the neural network using stochastic gradient descent methods, we leverage the propagation of chaos, which implies that particles
become asymptotically independent. Specifically, we train the neural network using a single simulation with large $N$, making our algorithm closer to 
the gradient-based training for a single-particle problem.
These extended algorithmic frameworks 
allow for a more general \emph{feedback functional} that depends 
not only on time and the particle's state, but also
\emph{on the empirical measure of the entire system}.  In principle, this expanded information set enables better generalization capabilities and more precise control strategies for problems in which the aggregate distribution of the population plays a critical role. 

Other algorithms for distribution-dependent FBSDEs have also been proposed  \cite{HHL}.  While there is a close connection between the FBSDE formulations of mean field control and their associated dynamic programming equations, these methods focus primarily on numerically solving the coupled systems of FBSDEs and are thus not directly comparable to our approach.  

\medskip
\noindent\emph{Theoretical results and convergence.}
From a theoretical viewpoint, a classical line of research (starting with Kac \cite{Kac1956} and subsequently developed in \cite{Funaki1984,McKean1966,Sznitman1991}) establishes \emph{propagation of chaos} for interacting particle systems, guaranteeing that the $N$-particle interaction system converges to a mean field limit as $N$ tends to infinity. In the context of mean field control, Lacker \cite{lacker} has provided general convergence results for $N$-particle problems to their corresponding MFC limits.  While the results in \cite{lacker} apply quite broadly, they do not directly yield the specific structure we exploit here.  Instead, in the particular setting of this paper, in Proposition \ref{viscosity}, we provide a simpler argument showing convergence of our single simulation approximation to the solution of the mean field control problem. 
The proof given in the Appendix is motivated by the techniques from the viscosity theory.

Additionally, in many applications, the mean field control problem is formulated in continuous time.  A standard time-discretization argument based on viscosity-solution theory \cite{FS}  can be used to justify the convergence of the discrete-time approximation to the continuous-time model.  
In this context one employs the classical methods developed by Barles-Souganidis \cite{BS} and Barles-Perthame \cite{BP} together with the comparison results of \cite{SY,SY1}.  We shall not reproduce all those classical steps, but we emphasize that our setting naturally fits into these well-established frameworks.  

\medskip
\noindent{\emph{Contributions of this paper.} 
We analyze a ML-based algorithm for mean field optimal control that:
\begin{itemize}
\item uses \emph{one large-population simulation} to generate training data, rather than the stochastic gradient approach which uses multiple simulations of smaller $N$ systems repeated in each epoch.  
\item trains a neural network policy that depends on the aggregate state distribution (via the empirical measure) as well as the individual player's state and time. We apply here flexible feature choice suggested by our universal approximation theorem.
\item achieves reliable convergence guarantees to the exact mean field control solution with a complete proof under general assumptions. In particular, we apply here the global character of our universal approximation theorem on Wasserstein spaces.
\item  tests on several representative MFC examples including a linear-quadratic model with an explicit solution, a Kuramoto mean field control problem exhibiting phase-transition behavior, a systematic-risk model for which comparison with existing deep-learning benchmarks and theoretical solution is available, and a crowd-aversion planning problem. These examples illustrate both the accuracy of the method and the computational efficiency gained by using one large-population simulation.
\end{itemize}
Besides offering computational efficiency, especially in settings where simulating a large population only once is practical, our method also demonstrates enhanced flexibility in capturing feedback policies sensitive to the distribution of the entire population.  The convergence proof, which we provide, adapts classical propagation-of-chaos arguments to our framework. We also use Rademacher complexity to bound the randomness in our convergence proof.}

\medskip
\noindent\emph{Organization of the paper.}
The rest of this paper is organized as follows. We summarize our notations in the rest of this section. In Section~\ref{sec:problem}, we present the mathematical formulation of the mean field control problem and the associated $N$-particle approximation.  Section~\ref{sec:algorithm} describes our numerical algorithm, including the neural network architecture for the feedback policy and the training procedure. A universal approximation theorem on Wasserstein spaces is  developed in Section~\ref{sec:UAT}. Numerical experiments illustrating the performance gains and efficacy of the proposed method are given in Section~\ref{sec:experiments}. In Section~\ref{sec:analysis}, we provide the convergence analysis of our numerical algorithm.  Concluding remarks and possible extensions are provided in Section~\ref{sec:conclusion}.

\subsection{Notations}
Let $E$ be a $d$-dimensional Euclidean space or a torus,  
$\cP(E)$ the set of all Borel probability measures on $E$. For $p \ge 1$, $\cP_p(E)$ is the set of probability measures 
with finite $p$-moments, and $\cP_\infty(E)$ consists of those
with compact support.
$\cW_p$ denotes the $p$-Wasserstein distance on $\cP_q(E)$ for $q\geq p$.
Fix $T\in \N$, and set 

$$
\cT:=\{0,\ldots,T\},\qquad  \text{and}
\qquad
\cX:=\cT\times E\times \cP_2(E).
$$
We endow $\cX$ with the product topology of the given topology on $E$ and the Wasserstein-2 metric topology.
Let $(\Omega, \mathbb F=(\cF_j)_{j=0,\ldots,T}, \mathbb P)$ be
a filtered probability space.  $\R^d$ denotes the $d$-dimensional Euclidean space
and $\cM(d,m)$ the set of all $d$ by $m$ matrices.
The noise process
$$
\boldsymbol{\xi}:=(\xi^{(1)},\xi^{(2)},\ldots)
\quad \text{with}\quad
\boldsymbol{\xi}^N:=(\xi^{(1)},\xi^{(2)},\ldots,\xi^{(N)})
$$
is an infinite i.i.d.~random sequence, where for every $i\in \N$,
the $i$-th random process $\xi^{(i)}:=(\xi^{(i)}_0,...,\xi^{(i)}_{T-1})$
with $\xi^{(i)}_j \in \R^m$
has independent components sampled from a 
fixed distribution $\gamma \in \cP(\R^m)$ and is adapted to the filtration $\F$. 
In many applications,
the common distribution of $\xi^{(i)}_j$ is Gaussian with zero mean
and small variance, representing the time discretization of 
Brownian motion.

For  $x=(x^{(1)},...,x^{(N)})\in E^N$, 
the corresponding empirical measure on $\cP(E)$ 
is given by 
$$
\mu^N(x):=\frac{1}{N}\sum_{i=1}^N \delta_{x^{(i)}},
$$
where $\delta_{x}$ is the Dirac measure located at $x$. For a random variable $X$, $\cL(X)$ denotes its distribution.

For any topological spaces $\cY,\cZ$, $\cC(\cY \to \cZ)$ is the set of all continuous functions 
$f:\cY \to \cZ$, $\cC(\cY):= \cC(\cY \to \R)$, and $\cC_b(\cY \to \cZ)$, $\cC_b(\cY)$ are the bounded ones. The 
\emph{control set} $A$ is a closed Euclidean space, and the set of \emph{feedback controls} $\cC$
is the collection of all functions  $\alpha \in \cC_b(\cX\rightarrow A)$
such that there is a constant $c(\alpha)$ depending on $\alpha$ and
for every   $(j,x,\mu), (j,y,\nu) \in \cX$,
\be
\label{eq:controlset}
|\alpha(j,x,\mu)| \le c(\alpha),\qquad
|\alpha(j,x,\mu)-\alpha(j,y,\nu)| \le c(\alpha) (|x-y|+\cW_1(\mu, \nu)).
\ee
For every function $\psi: \cX\times A\mapsto \R$ and a feedback control  $\alpha\in \cC$, set  
$$
\psi^\alpha(j,x,\mu):=\psi(j,x,\mu,\alpha(j,x,\mu)),
\qquad \forall\,(j,x,\mu)\in \cX.
$$
In the optimization literature 
$\psi^\alpha(\cdot)=\psi[\alpha](\cdot)$ is 
sometimes called the Nemytskii operator.

\section{Control problems}
\label{sec:problem}
In this section, we define the 
mean field control problem
we study and a particle system used 
to approximate it.  We assume that the following functions
are given and fixed throughout the paper:
$$
b: \cX \times A \mapsto E,\quad
\sigma :\cX \times A \mapsto \cM(d,m),\quad
\ell :\cX\times  A\mapsto \R,\quad
G : E \times \cP_2(E) \mapsto \R.
$$

\subsection{Mean field optimal control problem}
\label{sec:prob_formulation}

For a given feedback control 
$\alpha \in \cC$ 
and an initial distribution $\mu \in \cP_2(E)$, suppressing
the dependence on $\alpha$, the representative particle  
follows the following mean field equation
\begin{align}
   \label{eq:mfdynamics} 
   X_{j+1}&=X_j+b^\alpha(j,X_j,\cL(X_j))+\sigma^\alpha(j,X_j,\cL(X_j))\ \xi_j, 
    \quad j=0,\dots,T-1,
\end{align}
with an initial condition $X_0$ satisfying $\cL(X_0)=\mu$.
The solution $\boldsymbol{X}:=\boldsymbol{X}^\alpha=(X_0,\ldots,X_T)$
is a deterministic, albeit highly
complex, function of the random {{process $\xi$}} and the initial condition $X_0$. 
It is clear that this dependence is adapted. 
When we wish to emphasize this dependence we write 
$X_j({\xi},X_0)$. 

The path-wise cost function that we wish to minimize is defined by
$$
 \mathfrakJ(X_0,\alpha, {{\xi}}):=
 \sum_{j=0}^{T-1} \  \ell^\alpha(j, X_j({\xi},X_0),\cL(X_j))
+G(X_T({\xi},X_0), \cL(X_T)).
$$
Then, the mean field optimal control problem is 
to minimize the following cost function 
\be
\label{eq:jequal}
J(\mu, \alpha):=\E[ \mathfrakJ(X_0,\alpha, {{\xi}})]
= \E[\ \sum_{j=0}^{T-1} \ell^\alpha(j, X_j({\xi},X_0),\cL(X_j))
+G(X_T({\xi},X_0), \cL(X_T)\ ] ,
\ee
and the value function is given by
$$
v(\mu):=\inf_{\alpha\in \mathcal{C}} J(\mu, \alpha).
$$

\subsection{$N$-Particle System}
\label{sec:Nparticle}
For a fixed integer $N\in \mathbb N$,
we approximate the mean-field problem by an 
$N$-particle stochastic optimal control problem. 
We first fix the initial condition
$x^N=(x^{(1)},...,x^{(N)})\in E^N$.
Then, for a given feedback control $\alpha \in \cC$,
suppressing dependence on $\alpha$, we define
$$
\boldsymbol{X}^N:=(\boldsymbol{X}^{(1)},...,\boldsymbol{X}^{(N)}),
\qquad \text{with} \qquad
\boldsymbol{X}^{(i)}:=(X_0^{(i)},...,X_T^{(i)}),
\quad i=1,\ldots N,
$$
through the recursion with $X^N_0=x^N$,
\begin{equation}\label{eq:N-dynamics}
    X_{j+1}^{(i)}=X_j^{(i)}+b^\alpha(j,X_j^{(i)},\mu^N_j)
    +\sigma^\alpha(j,X_j^{(i)},\mu^N_j)\ \xi_j^{(i)},\quad
    j=0,...,T-1, \, i=1,\ldots,N, 
\end{equation}
where 
$$
X^N_j:=(X^{(1)}_j,\dots,X^{(N)}_j),\qquad \text{and}
\qquad
\mu^N_j := \mu^N(X^N_j).
$$
{{While the interaction in \eqref{eq:mfdynamics}
is through the law of $X_j$,
in \eqref{eq:N-dynamics} {the} empirical measure $\mu^N_j$ is used.
Therefore, $X_j^{(i)}$ depends on all components of the 
random sequence $\bsx=(\xi^{(1)}, \xi^{(2)},\ldots)$ and not only on $\xi^{(i)}$.}}  The $N$-particle optimization problem is to minimize the following cost function 
$$
J^N(x^N,\alpha):=\mathbb E\left[\frac{1}{N}\sum_{i=1}^N \sum_{j=0}^{T-1} 
\ell^\alpha(j,X^{(i)}_j,\mu^N_j)+G(\mu^N_T)\right],
$$
where the expectation is over the random sequences 
$\boldsymbol{\xi}$.
Then, the value function is given by
$$
v^N(x^N):=\inf_{\alpha\in \mathcal{C}} J^N(x^N,\alpha).
$$ 
\begin{remark}
\label{rem:motivate}
    The optimization problem is invariant under every permutation of indexes $i\in \{1,\dots,N\}$, 
    so the cost function $J^N$ and the value function $v^N$ are actually  functions of $\mu^N_0=\mu^N(x^N)$. 
    Moreover, the above particle system 
    is a natural approximation of the original problem.
    Indeed, by law of large numbers, 
    $$
    \E[\ell^{\alpha}(j,X,\cL(X_j))] \approx
    \E_{\mu^N(X^N_j)}[\ell^{\alpha}(j,X,\cL(X_j))]
    =\frac{1}{N}\sum_{i=1}^N \ell^\alpha(j,X^{(i)}_j,\cL(X_j)).
    $$
    We then choose the components of the random initial
    data $X_0^N$ independently from $\mu$, so that
    the propagation of chaos results apply to the 
    particle system \eqref{eq:N-dynamics}.  Consequently, 
    we may approximate $\cL(X_j)$ by the empirical 
    distribution $\mu^N_j$ so that  $J^N(X_0^N,\alpha)$ formally 
    approximates $J(\mu,\alpha)$.
    We thus expect $v^N(X_0^N)$ to be close
    to the main value function $v(\mu)$
    for such random sequences. 
    This formal argument
    motivates the definition of the $N$-particle problem
    and is rigorously proved in Section \ref{sec:analysis} and Appendix \ref{viscosity}.
\end{remark}

In our algorithm, as argued in the above remark,
we use random sequences $X_0^N$ to approximate $\mu$ and
additionally we fix a  random sequence $\boldsymbol{\xi}$. 
As in the mean-field case, $\boldsymbol{X}^N$ is also
a deterministic function of the random sequence 
$\bsx$ and the initial condition $x^N$. When we wish to emphasize this dependence we write 
$X^N_j(\boldsymbol{\xi},x^N)$, and $\mu^N_j(\boldsymbol{\xi},x^N)$. 
It is also clear that this dependence is adapted, i.e.,
$X^{N}_j(\boldsymbol{\xi},x^N)$ depends only on $\boldsymbol{\xi}_k$ for $k \le j-1$.
We introduce the  resulting 
(random) value functions of $\bsx$, $X_0^N$ that are
used centrally in our algorithm,
\begin{align}
\label{eq:jn}
\mathfrakJ^N(X^N_0,\alpha;\boldsymbol{\xi})&:= \frac{1}{N}\ \sum_{i=1}^N\sum_{j=0}^{T-1} 
 \ell^{\alpha}(j,X^{(i)}_j(\boldsymbol{\xi},X^N_0), \mu^N_j(\boldsymbol{\xi},X^N_0))+G(\mu^N_T(\boldsymbol{\xi},X^N_0)),\\
 \nonumber
 V^N(X^N_0;\boldsymbol{\xi})&:=\inf_{\alpha\in \mathcal{C}} 
\mathfrakJ^N(X^N_0,\alpha;\boldsymbol{\xi}).
\end{align}
We note that
$$
J^N(x^N,\alpha)=\E[ \mathfrakJ^N(X_0^N,\alpha; \bsx)\ |\ X_0^N=x^N].
$$
However, $\E[V^N(X_0^N;\xi)\mid X_0^N=x^N]$ is not equal
to $v^N(x^N)$.  Still in our analysis we show that they
are asymptotically close to each other in the sense stated
in Theorem~\ref{convergence} and Appendix~\ref{viscosity} below for appropriately
chosen sequences $x^N$ or $X_0^N$.

{\begin{remark}[Attainment of minimizers]
    The value functions $v(\mu)$, $v^N(x^N)$, and $V^N(X_0^N;\xi)$ introduced above are defined as infima over the admissible feedback class $\mathcal C$. Under the standing assumptions of this paper we do not claim that these infima are attained. Indeed, the class $\mathcal C$ is neither {$\sigma$-compact nor locally compact} in a topology under which the controlled dynamics and cost functional are lower semicontinuous; furthermore, the action space $A$ is only assumed to be a closed Euclidean space. Existence of optimal controls for stochastic control and McKean--Vlasov control problems typically requires additional compactness, coercivity, convexity, or relaxed-control assumptions; see, for instance, \cite{FS,CD,Bahlali_2017}. Thus, in the unrestricted problem over $\mathcal C$, optimality should be understood in the sense of approximate optimality.

    By contrast, the numerical optimization is performed after restricting the controls to a fixed hypothesis class $\mathcal N_K$, introduced in Section~\ref{sec:hypothesisclass}. For each fixed $K$, the parameter set $\Theta_K$ is compact and finite-dimensional, and the map $\theta\mapsto \alpha(\cdot;\theta)$ is assumed to be continuous. Under the continuity assumptions on the coefficients and costs, the maps
$$
\theta\mapsto J(\mu,\alpha(\cdot;\theta)),\qquad
\theta\mapsto J^N(x^N,\alpha(\cdot;\theta)),
$$
and, for fixed $(X_0^N,\xi)$,
$$
\theta\mapsto J^N(X_0^N,\alpha(\cdot;\theta);\xi)
$$
are continuous. Hence the restricted infima defining $v_K(\mu)$, $v_K^N(x^N)$, and $V_K^N(X_0^N,\xi)$ are attained. The convergence theorem below should therefore be interpreted as saying that, as $K$ and $N$ increase, these attainable minima over the hypothesis classes approximate the original mean-field value $v(\mu)$. We do not assert the existence of an optimizer over $\mathcal C$, nor convergence of minimizers in $\mathcal C$.
\end{remark}}

\section{Numerical algorithm}
\label{sec:algorithm}
In this section, we outline a neural network-based numerical algorithm to solve the mean field optimal control problem that is similar to the one studied in \cite{DLZ}. 
In view of Remark~\ref{rem:motivate} above, 
we fix a large $N$ and  simulate an
initial position $X^N_0$ from the distribution $\mu$
and the random sequences $\boldsymbol{\xi}^N$.
We then compute
the random value function $V^N(X^N_0;\boldsymbol{\xi})$
as an efficient approximation of $v(\mu)$.
One final ingredient
is the approximation of the optimal feedback
controls by neural networks so as to
efficiently compute the minimizers.
Section~\ref{sec:analysis} below,
justifies this approach 
by proving several convergence results, 
in particular, the convergence of
$V^N(X^N_0;\boldsymbol{\xi})$
to $v(\mu)$. 

\subsection{Hypothesis Class}
\label{sec:hypothesisclass}
In  our algorithm, we approximate the feedback controls
{{by a sequence of hypothesis classes.
In our numerical experiments, 
the hypothesis class is a set of neural networks.}}  We first introduce this construction more abstractly.  For $K\in \N$, let 
$$
\cN_K:=\{\alpha(\cdot;\theta):\mathcal{T} \times E\times \R^{m(K)} \rightarrow A \, | \, \theta\in \Theta_K\},
$$
be a sequence of hypothesis classes indexed by $K$,
where $\Theta_K$ is a finite-dimensional compact
parameter set and $m(K)$ is an increasing sequence
of positive integers. We assume that for each $K$,
elements of $\cN_K$ are uniformly bounded and Lipschitz
{{ with a uniform Lipschitz constant $K$}}.
To use any probability measure 
$\mu \in \cP(E)$ as an input to an element of the hypothesis
class, we need to approximate it by a finite-dimensional vector.
So, we further let 
$\cG:=\{g_1,g_2,\,...\} \subset \cC_b(E)$ be a separating class for $\cP_2(E)$, and for each $l$ set
$$
\cG_l(\mu):= \mu(g_l):=\int_E g_l(x)\,\mu(\d x),
\qquad \mu \in \cP(E).
$$
We then approximate $\mu$, by the finite-dimensional feature vector 
$(\cG_1(\mu),\ldots, \cG_{m(K)}(\mu))$
and  with an abuse of the notation, for any $(t,x,\mu)\in \cX$, and $\alpha(\cdot;\theta)\in \cN_K$, we set
$$
\alpha(t,x,\mu; \theta):=\alpha(t,x,(\cG_1(\mu),\ldots,
\cG_{m(K)}(\mu)); \theta).
$$
We first restrict the feedback controls to the
hypothesis class $\cN_K$ and set
\begin{equation}
\label{ee:vK}
v_K(\mu):= \inf_{\alpha \in \cN_K} J(\mu,\alpha),
\ \ \mu \in \cP(E),
\qquad
v_K^N(x^N):= \inf_{\alpha \in \cN_K} J^N(x^N,\alpha),
\ \ x^N \in E^N.
\end{equation}

{Since $\Theta_K$ is compact and the parameter-to-control map is continuous, the infima in the restricted problems over $\cN_K$ are attained whenever the corresponding objective function is continuous in $\theta$. In the sequel, these restricted minimizers are the objects computed by the numerical algorithm.}

To approximate the above problem
by an $N$-particle system, we simulate an
initial condition $X^N_0$ from the distribution $\mu$ and random sequences
$\boldsymbol{\xi}^N$.
By restricting to the feedback control set $\cN_K$, we define the (random) value function over the hypothesis class set by
$$
V^N_K(X_0^N,\boldsymbol{\xi}):=\inf_{\alpha\in \cN_K} 
\mathfrak J^N(X_0^N,\alpha;\boldsymbol{\xi}).
$$
Since  $\cN_K$ is a parameterized family of functions,
we introduce  the \emph{loss function} 
\begin{equation}
\label{eq:loss}
 L(\theta,X_0^N,\boldsymbol{\xi}):= \mathfrak J^N(X_0^N,\alpha(\cdot;\theta);\boldsymbol{\xi}),  
 \qquad \theta \in \Theta_K,
\end{equation}
and use the gradient descent based algorithm 
to minimize the loss function and train 
for an optimal parameter $\theta^*_K \in \Theta_K$. 
Then, in view of the universal approximation,
a uniform law of large numbers,
and propagation of chaos, when $N$ and $K$ are large enough, we formally  have
$$
L(\theta_K^*,X_0^N, \boldsymbol{\xi})
= V^N_K(X_0^N,\boldsymbol{\xi}) 
\approx v^N_K(X_0^N) \approx v_K(\mu)
\approx v(\mu).
$$
In Section~\ref{sec:analysis} below, we provide rigorous proofs
of the above approximations.

The main condition we place
on the hypothesis class $\cN_K$ is the classical universal approximation property.  
Recall that a feedback control  $\alpha \in \cC$ is a bounded continuous
function $\alpha :\cX \mapsto A$.

\begin{assumption}[Universal approximation]
\label{NN:UAT} 
     We assume that for every continuous function $\alpha\in \mathcal{C}$ and every $\epsilon>0$, for every compact set 
     $\mathfrak{C} \subset \cX$, there exists 
     $K= K(\eps,\mathfrak{C}) \in \N_+$ and $\alpha_K(\cdot;\theta)\in \cN_K$, such that 
    $$
    |\alpha(j,x,\mu;\theta)-\alpha_K(j,x,\mu)|\leq \epsilon,\quad \text{ for all } \,(j,x,\mu)\in \mathfrak{C}.
    $$
\end{assumption}

This local assumption is satisfied by many hypothesis classes. In particular, 
we provide a discussion of it for the class of  neural networks, 
and prove a global universal approximation  Theorem~\ref{UAT:Brhospace}  
in Section~\ref{sec:UAT} below.  Then,
the following result follows directly from these discussions and results. 

{
We define the weight function of quadratic growth by
    $$
    \rho(x) := 1 + \|x\|^2, \qquad x \in E,
    $$
and for $\mu \in \mathcal{P}_2(E)$, set
    $$
    R(\mu) := \int_E \rho(x)\,\mu(\d x).
    $$
Then, $R$ is a weight function on $\mathcal{P}_2(E)$ equipped with the weak topology. 
Finally we define  a weight on $E \times \mathcal{P}_2(E)$ by
    $$
    (\rho \vee R)(x,\mu) := \max\bigl(\rho(x), R(\mu)\bigr), \qquad (x,\mu)\in E\times\mathcal{P}_2(E).
    $$
{For a positive weight $w$ on a topological space $S$, i.e.~a function $w$ where inverse images of intervals are compact in $S$, let
$$
\mathcal B^w(S)
:=
\overline{C_b(S)}^{\|\cdot\|_w},
\qquad
\|f\|_w
:=
\sup_{z\in S}\frac{|f(z)|}{w(z)}.
$$
In the present setting, $\cP_2(E)$ is equipped with the weak topology, and we set
$$
\mathcal B^{\rho\vee R}(\mathcal X)
:=
\left\{
\alpha:\mathcal T\times E\times \cP_2(E)\to A:
\max_{j\in\mathcal T}
\sup_{(x,\mu)\in E\times P_2(E)}
\frac{|\alpha(j,x,\mu)|}{(\rho\vee R)(x,\mu)}
<\infty
\right\},
$$
where for each $j\in \cT$,
$(x,\mu)\mapsto \alpha(j,x,\mu)$ belongs to
$\mathcal B^{\rho\vee R}(E\times P_2(E))$. A more detailed definition of $\cB^{\rho \vee R}$
is given in Section~\ref{sec:UAT} below.}

\begin{theorem}
\label{thm:ourUAT}
For every $\epsilon>0$ and every
$\alpha \in \cB^{\rho \vee R}\bigl( E \times \mathcal{P}_2(E)\bigr)$,
there exists $K = K(\epsilon) \in \mathbb{N}_+$ and $\alpha_K(\cdot;\theta) \in \cN_K$ such that
$$
\bigl|\alpha(j,x,\mu) - \alpha_K(j,x,\mu;\theta)\bigr|
\;\leq\;
\epsilon \,\max\Bigl(\rho(x), \int_E \rho(y)\,\mu(\d y)\Bigr),
\qquad \text{for all } (j,x,\mu)\in \cT \times E \times \mathcal{P}_2(E).
$$
\end{theorem}}
We emphasize the global character of the above assertion
and the slightly weaker topology on $\mathcal{P}(E)$ in the assumptions on $\alpha$.

\subsection{Algorithm}

We fix large integers $K$, $N$, $m(K)$, and choose a separating class of functions $\cG:=\{g_1,g_2,...\}$ 
based on the structure of the problem. We  use the $m(K)$-dimensional vector $(\cG_1(\mu),\ldots,\cG_{m(K)}(\mu))\in\R^{m(K)}$ to approximate any $\mu \in \cP(E)$. 
Let $\cN_K$ be the set of neural networks induced from this approximation
with a compact parameter set $\Theta_K$. Then,
$$
\cN_K:=\{\alpha(\cdot;\theta):\mathcal{T} \times E\times \R^{m(K)}
\rightarrow A \,| \;\theta\in \Theta_K\}
$$
is a set of neural networks which are $K$-Lipschitz continuous and bounded by $K$. 
Then, any element $\alpha(\cdot,\theta)\in\cN_K$ is a 
feedforward neural network with input dimension $1+d+m(K)$, weight parameters $\theta\in \Theta_K$, where $\Theta_K\subset {\mathbb R}^{d(K)}$ 
is the parameter space with dimension $d(K)$. As $K$ goes to infinity, $d(K)$  also converges to infinity.

We first approximate the mean field dynamics with an $N$-particle system defined in Section 
\ref{sec:Nparticle}. 
That is, for a given distribution $\mu\in \cP_2(E)$, 
we simulate initial random variables
$(X_0^{(1)},..., X_0^{(N)})$ independently and identically distributed according to $\mu$,  then 
additionally generate independent random sequences
$\boldsymbol{\xi}^N=(\boldsymbol{\xi}^{(1)},\dots,\boldsymbol{\xi}^{(N)})$ from $\gamma$.
We use the loss function $L(\theta_K,X^N_0,\boldsymbol{\xi})$ of \eqref{eq:loss}, and then  train the parameters with gradient descent based algorithm.  
The pseudo-algorithm is given as follows.
\vspace{5pt}

\begin{algorithm}[h!]
\caption{Learning Algorithm for MFC}\label{alg:mfc}
\begin{algorithmic}
    \State Fix $N$, neural network $\cN_K$, initial condition $\mu \in \cP(E)$,
    and a noise distribution $\gamma \in \cP(\R^m)$;
    \State {\bf{Input\ data:}}  A neural network $\alpha(\cdot;\theta)\in \cN_K$;
    \State ${\bf{Initialization:}}$ Simulate random sequences $\boldsymbol{\xi}^N$
    and i.i.d.~samples $X_0^N=(X_0^{(1)},...,X_0^{(N)})$ from $\mu$;
    \For{n=1,2,...}
    \For{$j \gets 1$ to $T$}
        \For{$i\gets 1$ to $N$}
            
            \State Compute $\{\cG_1(\mu^N_{j-1}),\ldots,\cG_{m(K)}(\mu^N_{j-1})\}$ as in Remark \ref{remark:fourier} below;
            \State Update $X^{(i)}_j$ as in equation \eqref{eq:N-dynamics};
        \EndFor
    \EndFor
    \State Compute the loss function $L(\theta,X^N_0,\boldsymbol{\xi})$;
    \State Compute $d:=\nabla_\theta L(\theta,X_0^N,\boldsymbol{\xi}) $ as in \eqref{eq:loss};
    \State Update $\theta$ with $d$ ;
    \EndFor
    \State Return $\theta^*=\theta$.
    \State $\alpha(\cdot,\theta^*)$  is an approximation of the 
    optimal control;
    \State 
    $L(\theta^*,X_0^N,\boldsymbol{\xi})$ is an approximation of 
    the value function $v(\mu)$.
\end{algorithmic}
\end{algorithm}

\begin{remark}\label{remark:fourier}
    Thanks to the special structure of the empirical distribution
    $\mu^N_j=\mu^N(X^N_j)$, at each time step $j\in\cT$, we can compute the integrations efficiently by
    $$
    \cG_l(\mu^N_j)=\int_E g_l(x)\mu^N_j(\d x) =\frac{1}{N}\sum_{i=1}^N g_l(X^{(i)}_j),\quad l=1,\ldots,m(K).
    $$
\end{remark}
\begin{remark}
\label{rmk: seprarating class}
    There is flexibility in choosing the separating class of functions $\cG$. For example, in the Kuramoto model 
    in Section~\ref{example: Kuramoto}, as the underlying space $E$ is the torus, we employ the
    Fourier basis. For the other two examples we
    use simple polynomials (on bounded domains).  
    \cite{PW_wasserstein} considers piecewise linear functions as the separating class of functions.
    In more complicated problems, one can exploit the learning ability of the neural networks, 
    and let the neural network learn the optimal mapping $\cG^{m(K)}:\cP_2(E)\rightarrow \R^{m(K)}$ through training. 
\end{remark}

In our numerical examples, we use neural networks with 2 to 3 hidden layers, each containing 32 to 256 neurons. The activation functions are either $\tanh$ or $\operatorname{ReLU}$, and training is performed using the ADAM gradient method with 6000 to 10000 iterations. All numerical results were obtained using a personal MacBook laptop equipped with a 16 GB Apple M1 Pro chip.

\section{Universal Approximation on $\cX$}
\label{sec:UAT}

We shall develop universal approximation on Wasserstein spaces here. Our presentation deviates from \cite{PW_wasserstein} and is of independent interest. It mainly relies on insights from \cite{CST:2025}.

First some general remarks: if $S$ is a Polish space, i.e.~a completely metrizable separable space, 
then the probability measures over $S$ with the topology which metrizes weak convergence of sequences (weak convergence in turn is the trace of the weak star topology on the dual of $\cC_b(S)$), denoted by $\mathcal{P}(S)$, 
form a Polish space, too. In other words: there are complete metrics metrizing weak convergence on Polish spaces, e.g.
$$
M_d(\mu,\nu):=\inf_{X \sim \mu, \, Y \sim \nu} E[d(X,Y)]
$$
for any bounded complete metric $d$ on $S$. Being {Polish} is inherited by probability measures on a {Polish} space, but not being locally compact on locally compact spaces. We show in the sequel a second inheritance property in case of weighted spaces. 

If $S$ is compact, then $P(S)$ is compact, too, and the situation is fairly simple. 
In this case $\cC_b(S)=\cC(S)$ is separable and the weak star topology is metrizable on the unit ball which contains $\mathcal{P}(S)$. 
Since $\cC(S)$ is separable, we can choose countable families of functions $(f_i)$, 
which are point separating. We shall always assume that they take values in $[0,1]$. 
The map $S \ni x \mapsto (f_i(x)) \in [0,1]^\infty$ is continuous, injective and 
by basic topology a homeomorphism onto its compact image. This means in particular that any separable family is actually convergence characterizing. This yields a numerical model for $S$, which allows to pull back results on polynomial approximation and on neural networks with input space $[0,1]^\infty$ to $S$. By the same reasoning we can construct a numerical model on $\mathcal{P}(S)$, namely $ \mu \mapsto (\int_S g_i(x) \mu(dx))$, where $(g_i)$ is a countable family of functions, which characterizes weak convergence of measures on the compact $S$. Again we can pull back results on polynomial approximation from $[0,1]^\infty$ to $\mathcal{P}(S)$.
\begin{example}
Let $S$ be compact and let $(f_i)$ be a point separating family on $S$. Let $y \mapsto \phi(y)$ be a neural network taking inputs finitely many $(x \mapsto f_i(x))$. Then the set of all neural networks of this type is dense in $C(S)$ and has a countable subset $(\phi_j)$ (each taking values in $[0,1]$), whose span is again dense. Take, for instance, the neural networks with rational weights and appropriate values on $S$. Then by the same arguments $ (\mu \mapsto \int_S \phi_j (x) \mu(\d x)) $ is point separating and convergence characterizing on $\mathcal{P}(S)$, and neural networks taking finitely many $(\mu \mapsto \int_S \phi_j (x) \mu(\d x) )$ as inputs are dense in $\cC(\mathcal{P}(S))$. In other words: for every continuous $ F $ on $\mathcal{P}(S) $ with real values and for every $\epsilon > 0 $ there is a neural network $ \psi(\langle \phi,\mu) \rangle) $ on $ \langle \phi, \mu \rangle := (\int_S \phi_j (x) \mu(\d x))_j $, such that $ | F(\mu) - \psi(\langle \phi, \mu \rangle)| < \epsilon $ for all $\mu \in \mathcal{P}(S)$.
\end{example}

When $S$ is not compact, the situation is more complicated. It starts with the fact that in this case 
$\cC_b(S)$ is \emph{not} separable and therefore the weak star topology is \emph{not} metrizable on its unit ball. It is already remarkable that even though the weak star topology, which determines weak convergence of measures, is not metrizable, there is a metric metrizing the weakly converging sequences of measures. Furthermore $S$ is compact if and only if $\mathcal{P}(S)$ is compact. If $S$ is not compact $\mathcal{P}(S)$ is neither locally compact nor $\sigma$-compact making it difficult to prove approximation theorems. This is simply due to the fact that a metric ball in a transport metric can contain measures whose supports are far apart, such that no weakly converging subsequences need to exist.

One {workaround} for this problem is to use a slightly weaker version of weak convergence, namely vague convergence, which works in case $S$ is locally compact: one applies the one point compactification of $S$ and the topology inherited from the compact space $ S \cup \{\infty\}$. This leads to loss of mass to $\infty$ and since $\mathcal{P}(S)$ is not locally compact adding one point $\delta_\infty$ does not compactify $\mathcal{P}(S)$ with its topology of weak convergence, since obviously $P(S)$ is not locally compact even if $S$ is locally compact but not compact.

A more canonical approach, which fits better approximation theoretic questions (and ultimately questions from machine learning), works with weight functions. Assume that $S$ is {Polish} with a weight function $\rho : S \to \mathbb{R}_{>0}$, {i.e.~inverse images of compact intervals under $\rho$ are compact in $S$}. Consider $\cB^\rho(S)$ the closure of $\cC_b(S)$ with respect to the weighted uniform norm $\| f \|_{\rho} := \sup_{x \in S} \frac{|f(x)|}{\rho(x)}$. This is a Banach space whose dual can be identified with Borel measures $\mu$ such that 
$ \int_S \rho(x) |\mu|(\d x) < \infty $, since any bounded linear functional on $\cB^\rho(S) $ 
can be identified with $ f \mapsto \int_S f(x) \mu(\d x) $ {defined} for some {unique} measure $\mu$. Consider now $\mathcal{P}_\rho(S)$ the set of probability measures $\mu \in \mathcal{P}(S) $ such that the $\rho$ moment exists, in other words the set of probability measures in the dual space of $\cB^\rho(S)$. The weak star topology on $\mathcal{P}_\rho(S) $ apparently has a weight function again, namely 
$$ 
R(\mu) := \mu \mapsto \int_S \rho(x) \mu (\d x) 
$$
since $|\mu|=\mu$ for positive measures and since the strong balls in the dual space are compact by Banach-Alaoglu. Whence we have inheritance. Notice also that probability measures with $\rho$ moment bounded by a constant $M$ are metrizable with this topology, since $\cB^\rho(S)$ is separable. We have the following proposition.

{\begin{proposition}[Compact weighted-moment sublevels]
\label{prop:weighted-measure-compactness}
For every $M>0$, the set
$$
  K_M
  :=
  \left\{\mu\in\mathcal P_\rho(S):R(\mu)\le M\right\}
$$
is compact in the relative weak star topology.  Moreover, this topology is
metrizable on $K_M$.
\end{proposition}

\begin{proof}
For a positive measure $\mu$,
$$
  R(\mu)
  =
  \|\mu\|_{\mathcal B_\rho(S)^*}.
$$
Hence $K_M$ is the intersection of $\mathcal P_\rho(S)$ with the
closed ball of radius $M$ in $\mathcal B_\rho(S)^*$.  Banach--Alaoglu
implies that this dual ball is compact in the weak star topology.

It remains to check that $\mathcal P_\rho(S)$ is weak star closed in
$\mathcal M_\rho(S)$.  Suppose that a net
$(\mu_\alpha)\subset\mathcal P_\rho(S)$ converges weak star to
$\lambda\in\mathcal M_\rho(S)$.  Since $1\in\mathcal B_\rho(S)$,
$
  \lambda(S)=\langle 1,\lambda\rangle
  =\lim_\alpha\langle 1,\mu_\alpha\rangle=1.
$
Moreover, for every nonnegative $g\in C_b(S)$,
$$
  \int_S g\,\d\lambda
  =
  \lim_\alpha\int_S g\,\d\mu_\alpha
  \ge 0.
$$
Thus $\lambda$ is a positive measure of total mass one, and therefore
$\lambda\in\mathcal P_\rho(S)$.  This proves compactness of $K_M$.
Since $\mathcal B_\rho(S)$ is separable, the weak star topology on a
norm-bounded subset of its dual is metrizable.
\end{proof}}

Therefore we can formulate an analogous Stone-Weierstrass and universal  approximation result as above in this setting. Let $(f_i)$ be a separating family on $S$ taking values in $[0,1]$, which is convergence characterizing on compacts. Let again $x \mapsto \phi(x)$ be a neural network on inputs 
$(x \mapsto f_i(x))$. Then the set of all neural networks of this type is dense in 
$\cB^\rho(S)$ and has a countable subset $(\phi_j)$ (each taking values in $[0,1]$), whose span is again dense. Take again the neural networks with rational weights and appropriate {value set $[0,1]$} on $S$. Then, by the same arguments 
$ (\mu \mapsto \int_S \phi_j (x) \mu(\d x)) $ is point separating and convergence characterizing on 
$\mathcal{P}_\rho(S)$ on compacts, and neural networks with inputs 
$(\mu \mapsto \int_S \phi_j (x) \mu(\d x) )_j$ are dense in $\cB^R(\mathcal{P}_\rho(S))$. In other words: 
\begin{theorem}\label{UAT:Brhospace}
For every continuous $ F \in \cB^R(\mathcal{P}_\rho(S)) $ with input space 
$\mathcal{P}_\rho(S) $ and for every $\epsilon > 0 $ there is a neural network 
$ \psi(\langle \phi,\mu) \rangle) $ {taking 
$ \langle \phi, \mu \rangle := (\int_S \phi_j (x) \mu(\d x))_j $ as inputs}, 
such that $ | F(\mu) - \psi(\langle \phi, \mu \rangle)| < \epsilon R(\mu) $ for all $\mu \in \mathcal{P}_\rho(S)$.
\end{theorem}

Notice that in case of $S=\mathbb{R}^d$ the set $\mathcal{P}_\rho(S)$ coincides for $\rho(x)= 1 + \|x\|^p$ with the Wasserstein space $\mathcal{W}_p(\mathbb{R}^d)$, however, with a slightly weaker metric on sets with bounded $p$-th moment, namely just weak convergence. Continuity in Theorem \ref{UAT:Brhospace} refers to this metric. Notice also that Borel and Baire $\sigma$-algebra with respect to these different metrics coincide.

\begin{remark}
In contrast to \cite{PW_wasserstein}, where a modulus of continuity is needed and neural networks depend on $K^*$ (even though this is not stated), we obtain a global weighted universal approximation theorem and a precise characterization of all functions which can be approximated.
\end{remark}

{\begin{remark}
    The {Universal Approximation Theorem~\ref{thm:ourUAT}} is stated qualitatively. Obtaining an explicit rate in $K$ requires additional regularity assumptions which are not imposed in the present paper. If the target feedback depends only on finitely many features and the resulting finite-dimensional function has Sobolev, H\"older, Barron, or piecewise smooth regularity, then some existing quantitative neural-network approximation estimates can be applied, e.g.,
\cite{yarotsky2017error,petersen2018optimal}.
\end{remark}}

Consider now a finite measure $\eta$ on $S$ with 
$R(\eta) = \int_S \rho(x) \eta(\d x) < \infty$, i.e.~finite $\rho$ moment. 
Then, $\cB^\rho(S)$ is dense in $L^1(\eta)$ and we obtain universal approximation 
with respect to $L^1$-norm. Of course an $L^p$ version of this result exists 
by considering the weight $\rho^{p}$ instead of $\rho$.

\begin{example}
Consider the {Polish} {space} $S:=\mathbb{R}^d \times \mathcal{W}_2(\mathbb{R}^d)$ 
and the measures $ \eta (\d y,\d \mu) := \mu(\d y) \otimes \nu(\d \mu)$ such that
$$
\int_S \|y\|^2 \mu(\d y) \nu(\d \mu) < \infty \, .
$$
Then neural networks of the type $ \mu \mapsto \psi( \langle \phi, \mu\rangle) $ are dense in $L^2(\eta)$, since the weighted functions with respect to the weight $ \mu \mapsto \sqrt{\int_{\mathbb{R}^d} \| x \|^2 \mu(dx)} $ are dense in $L^2(\eta)$.

Notice that this could also be seen by a martingale argument, where time is indexed by a countable index set coming from $(\phi_i)$.
\end{example}

\section{Numerical results}
\label{sec:experiments}
In this subsection, we test our learning algorithm on three mean field models originally formulated in continuous time. By applying an equidistant time discretization, we obtain corresponding discrete-time problems that converge to the original continuous-time problems as the discretization time step tends to zero.

\subsection{Linear quadratic problem}

{{We start with the following simple
linear quadratic problem with an explicit solution
as a benchmark.  The solution is a function 
of the first two moments of the distribution.  Although we do not encode this 
into our algorithm, it learns this finite-dimensional
structure quickly.}}  

In this setting, 
the underlying space $E=\R^d$, and the MFC problem with the interaction parameter $\kappa>0$ is
$$
v(\mu):=\inf _{\alpha \in \mathcal{A}} \mathbb{E} \int_0^{\infty} e^{-\beta t}\left(\frac{1}{2}\left|\alpha_t\right|^2+\kappa 
\operatorname{Var} \left(\mathcal{L}\left(X_t^\alpha\right)\right)\right) \mathrm{d} t,
$$
where $d X_t=\alpha_t d t+\sigma d W_t$, and $\cL(X_0)=\mu$. The dynamic programming equation is
$$
\beta v(\mu)=-\frac{1}{2} \mu\left(\left|\nabla_x \delta_\mu v(\mu)(\cdot)\right|^2\right)+\frac{\sigma^2}{2} \mu\left(\Delta_x \delta_\mu v(\mu)(\cdot)\right)+\kappa \operatorname{Var}(\mu).
$$
This is a linear quadratic formulation, and there is a unique analytic solution 
$$
v(\mu)=a\operatorname{Var}(\mu)+b,
$$
where
$$
a=\frac{\sqrt{8 \kappa+\beta^2}-\beta}{4}>0, \quad b=\frac{\sigma^2  a}{\beta}>0.
$$
The stationary minimizers are Gaussian distributions with arbitrary mean and variance $\frac{\sigma^2}{4a}$.

We apply our numerical algorithm to this linear quadratic problem, and use the analytic solution as our benchmark. In this experiment, we set $\kappa=1$, $\sigma=1$, $\beta=1$, then $a=0.5$, $b=0.5$, and the stationary distribution is $N(0,0.5)$. In the code, we use the ReLU activation function with 2 hidden layers, each with 32 neurons. We take simple polynomials as the separating class of functions, and use 10-dimensional feature vectors to approximate the measures. We use the Pytorch package and the Adam algorithm for stochastic optimization to train the parameters.

\begin{figure}[h]
    \begin{subfigure}[b]{0.5\textwidth} 
        \centering
        \includegraphics[width=\textwidth]{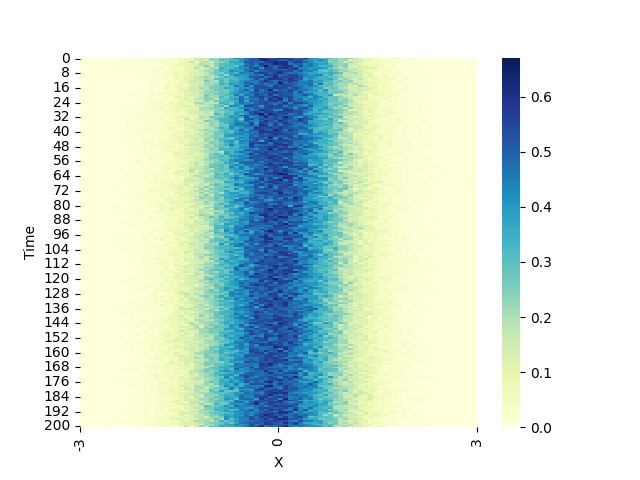}
        \caption{Stationary solution}
    \end{subfigure}
    \begin{subfigure}[b]{0.5\textwidth}
        \centering
        \includegraphics[width=\textwidth]{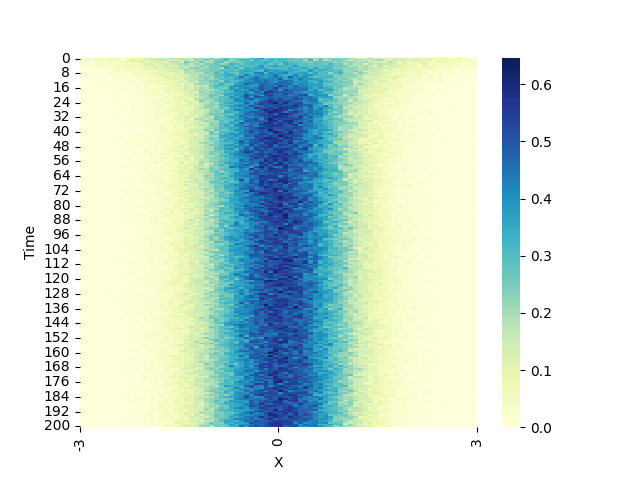}
        \caption{Dynamic solution}
    \end{subfigure}
    \caption{Linear Quadratic MFC}
    \label{fig:LQ-gaussian}
\end{figure}

 Figure \ref{fig:LQ-gaussian} shows two optimal flows with different initial distributions obtained from our numerical method. The left is the optimal flow starting with the stationary distribution $\mu=N(0,0.5)$. Our result shows that this optimal flow is constant in time. The right one is the optimal flow with initial distribution $\mu=N(0,9/4)$, and it converges to the stationary solution in time. The results align with previous analysis, showing that the Gaussian distribution with variance $0.5$ is the stationary solution.
\subsection{Kuramoto mean field control problem}
\label{example: Kuramoto}

{{Our second example is
related to the Kuramoto mean field game problem studied in \cite{CCS}
which exhibits
phase transition as the appropriate interaction parameter varies.
The problem we numerically study, the Kuramoto mean field control problem, 
is its potential and hence its minimizers
are Nash equilibria (NE) for the game, \cite{HS}.  
This connection allows us to test our numerical algorithm,
taking the analysis in \cite{CCS} as a benchmark.}}

{{We emphasize  that although the interaction that we define below
seems to have only first Fourier modes of the 
distribution $\mu$, the value function and consequently,
the optimal feedback control depend on
\emph{all Fourier modes} of $\mu$.  This makes the 
problem infinite dimensional without an explicit solution,
and it is approximated by
a finite but high-dimensional $N$-particle system.}}

The infinite-horizon, Kuramoto mean field control problem 
with state space $E=\T$ is defined as
$$
v(\mu):=\inf_{\alpha \in \cA} J(\alpha):=\inf_{\alpha\in \cA}\mathbb E \int_0^\infty e^{-\beta t}[\kappa \Phi(\mu_t)+\frac{1}{2}\alpha_t^2]\d t,
$$
where 
$$
\Phi(\mu):=-\frac{1}{2}\int_{{\mathbb T}^2} \cos{(x-y)} \mu(\d x)\mu(\d y),
$$
and $\mu_t=\text{Law}(X_t)$,  $\d X_t=\alpha_t \d t+\sigma \d W_t$, $X_0\thicksim \mu$. $\kappa>0$ is the interaction parameter.

The following phase transition phenomenon is proved in \cite{CCS}. 
\begin{theorem}
\label{thm:Kuramoto}
Let $\kappa_c:=\beta \sigma^2+\frac{1}{2}\sigma^4$.

\begin{enumerate}
    \item For all interaction parameters $\kappa>\kappa_c$, there are non-uniform stationary Nash equilibria;
    \item If $\kappa<\frac{1}{4}\beta\sigma^2$, then any Nash equilibria $\mu=(\mu_t)_{t\geq 0}$ with interaction parameter $\kappa$ converges in law to the uniform distribution, i.e., $\mu_t\rightarrow U$, as $t\rightarrow\infty$.
\end{enumerate}
\end{theorem}

The following numerical results are obtained from a fully-connected feedforward neural network. When $E=\T^d$, we take the Fourier basis as 
the separating class of functions, and with a slight change of notation, let the separating class be indexed by multi-indices and given by
$\{(\cos(2\pi l\cdot x), \sin(2\pi l\cdot x)) : l \in \Z^d\}$. We replace $\cG_l$ by
\begin{align*}
    C_l(\mu):=\int_E \cos{(2\pi l\cdot x)}\mu(\d x),\quad
    S_l(\mu):=\int_E \sin{(2\pi l\cdot x)}\mu(\d x),
    \quad l \in \Z^d,\ \mu \in \cP_2(E).
\end{align*}

In the code, we use the ReLU activation function with 2 hidden layers, each with  32 neurons, and take 
$m(K)$ to be 20. We use the Pytorch package and the Adam function for stochastic optimization to train the parameters.
The infinite time horizon is truncated at $T=20$, and then we discretize the time horizon in $T_n=400$ uniform segments. We fix $\beta=1$ and $\sigma=1$, thus $\kappa_c=1.5$ and take $N=3000$. The algorithm is tested with different $\kappa$ and different initial distributions. 

\subsubsection{Sub-critical case}
\begin{figure}[ht!]
    \begin{subfigure}[b]{0.5\textwidth} 
        \centering
        \includegraphics[width=\textwidth]{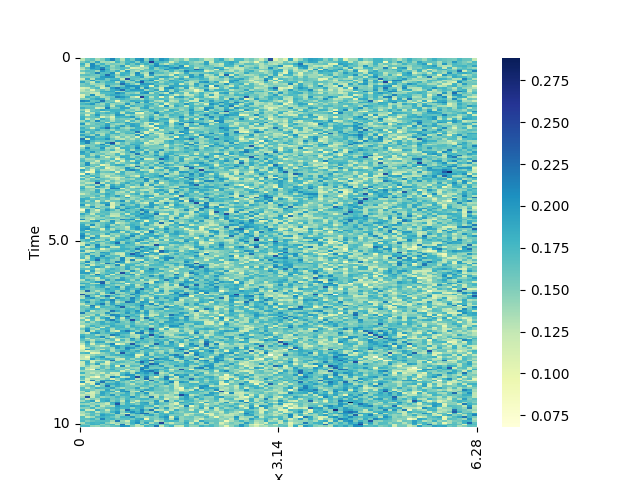}
        \caption{Heatmap of the minimizer}
        \label{fig:sub1}
    \end{subfigure}
    \begin{subfigure}[b]{0.5\textwidth}
        \centering
        \includegraphics[width=\textwidth]{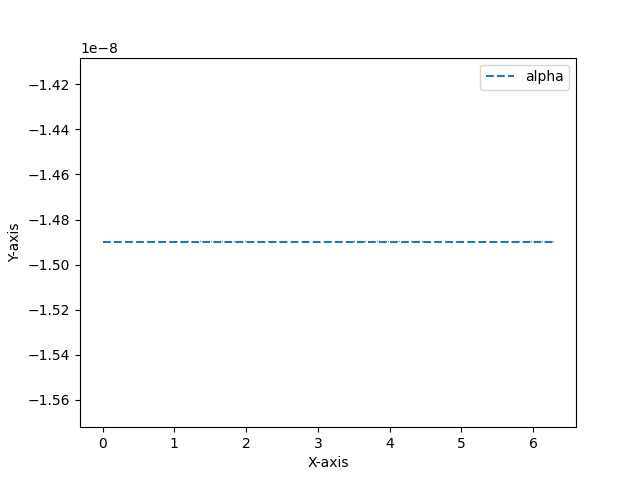}
        \caption{Optimal control at time T}
        \label{fig:sub2}
    \end{subfigure}
    \caption{Initial distribution: uniform distribution}
    \label{fig:uni_k0.8}
    
\end{figure}

\begin{figure}[h]
   \begin{subfigure}[b]{0.5\textwidth}
        \centering
        \includegraphics[width=\textwidth]{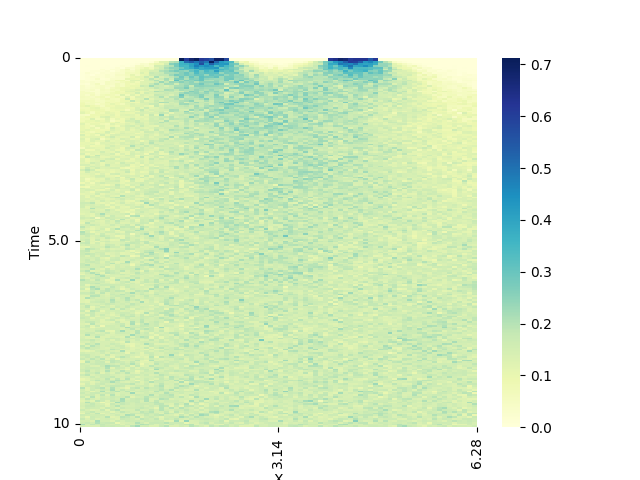}
        \caption{Heatmap of the minimizer}
    \end{subfigure}
    \begin{subfigure}[b]{0.5\textwidth}
        \centering
        \includegraphics[width=\textwidth]{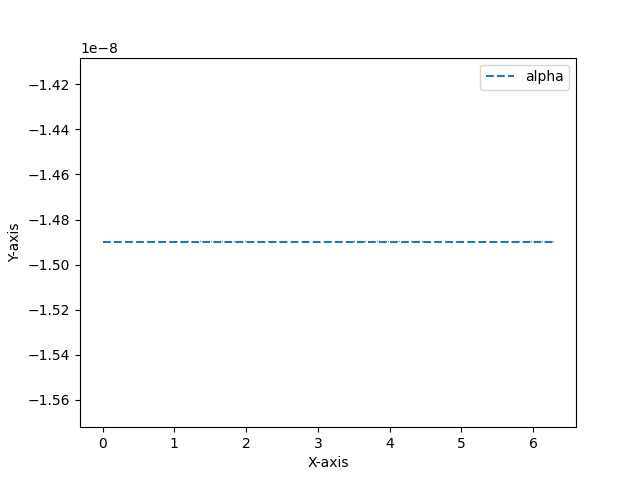}
        \caption{Optimal control at time T}
        \label{fig:sub4}
    \end{subfigure}
    \caption{Initial distribution: Two-cluster distribution}
    \label{fig:tc_k0.8}
\end{figure}

We first take $\kappa=0.2<\kappa_c$. Figure \ref{fig:uni_k0.8} and Figure \ref{fig:tc_k0.8} show the optimal flow in time and the optimal control at terminal time $T$. Figure \ref{fig:uni_k0.8} shows that when starting with uniform distribution, the optimal flow is a constant in time. Figure \ref{fig:tc_k0.8} shows the case with a two-cluster initial distribution, and the optimal flow converges to the uniform distribution in time. For both cases, the optimal control at the terminal time is approximately 0. The numerical results correspond with the sub-critical case in Theorem \ref{thm:Kuramoto}.

\subsubsection{Super-critical case}
\begin{figure}[h!]
    \begin{subfigure}[b]{0.5\textwidth} 
        \centering
        \includegraphics[width=\textwidth]{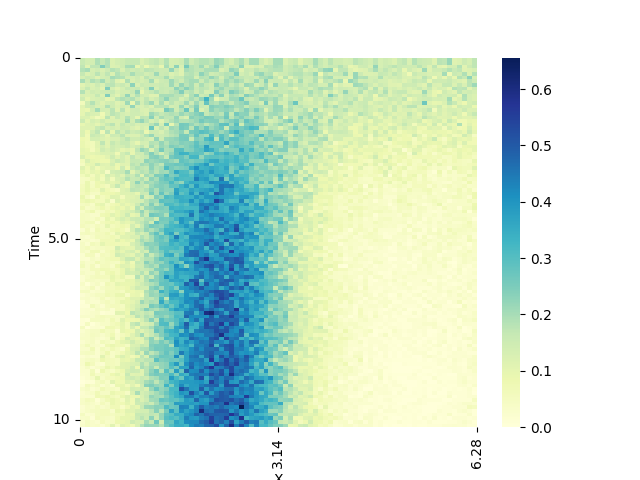}
        \caption{Heatmap of the minimizer}

        \label{fig:sub5}
    \end{subfigure}
    \begin{subfigure}[b]{0.5\textwidth}
        \centering
        \includegraphics[width=\textwidth]{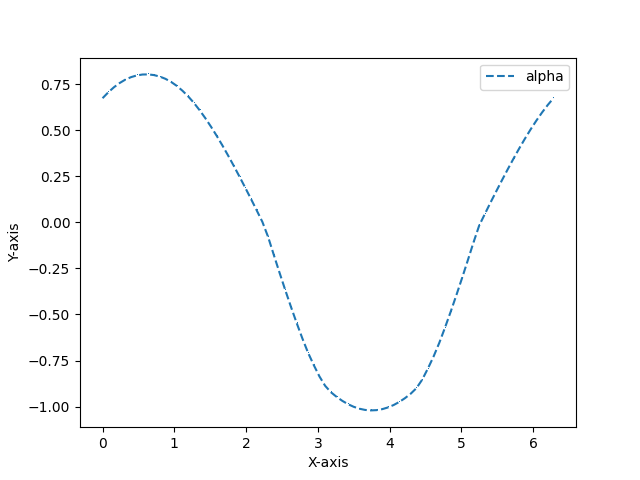}
        \caption{Optimal control at time T}
        \label{fig:sub6}
    \end{subfigure}
    \caption{Initial distribution: uniform distribution}
    \label{fig:uni_k2.5}
    
\end{figure}

\begin{figure}[h]
   \begin{subfigure}[b]{0.5\textwidth}
        \centering
        \includegraphics[width=\textwidth]{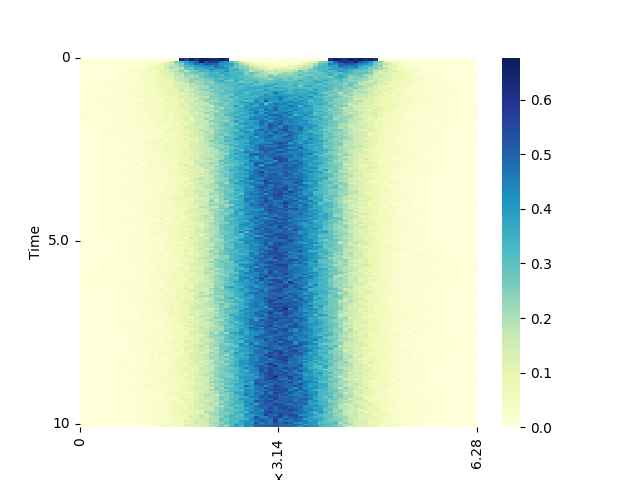}
        \caption{Heatmap of the minimizer}
    \end{subfigure}
    \begin{subfigure}[b]{0.5\textwidth}
        \centering
        \includegraphics[width=\textwidth]{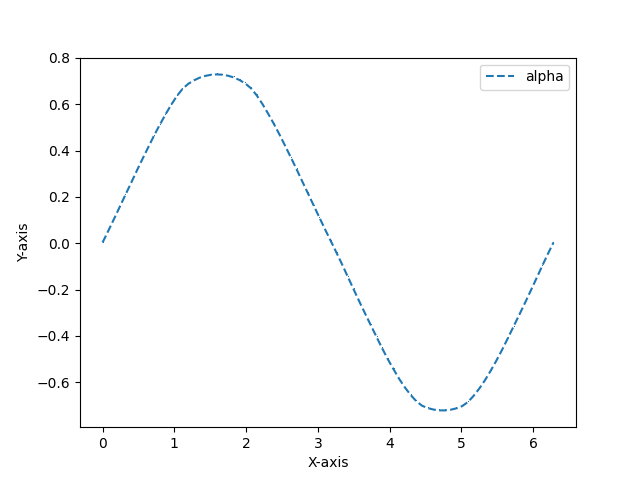}
        \caption{Optimal control at time T}
    \end{subfigure}
    \caption{Initial distribution: Two-cluster distribution}
    \label{fig:tc_k2.5}
\end{figure}

For the super-critical case, we take $\kappa=2.5>\kappa_c$.  Figure \ref{fig:uni_k2.5} and Figure \ref{fig:tc_k2.5} show the optimal flows and the optimal controls at terminal time $T$ with the two previous  initial distributions.
Notice that the optimal control at time $T$ is no longer 0, and although still a Nash equilibrium for the Kuramoto Mean field game, the uniform distribution is not a stationary minimizer for the Kuramoto mean field control problem. For both cases, synchronization happens and the optimal flows converge.

\begin{remark}
    This example also illustrates that for a potential mean field game, which is always numerically difficult and expensive to solve, 
    we can obtain one Nash Equilibrium by solving the corresponding mean field control problem with our numerical algorithm.
\end{remark}

\subsubsection{Generalization of the neural network}
\begin{figure}[h]
   \begin{subfigure}[b]{0.5\textwidth}
        \centering
        \includegraphics[width=\textwidth]{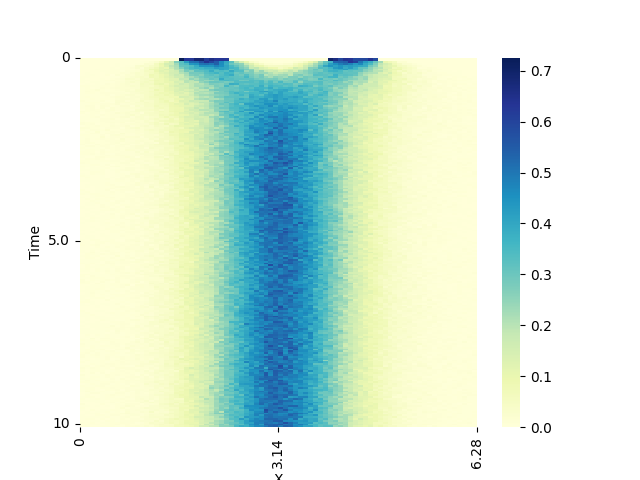}
        \caption{$\kappa=3$}
    \end{subfigure}
    \begin{subfigure}[b]{0.5\textwidth}
        \centering
        \includegraphics[width=\textwidth]{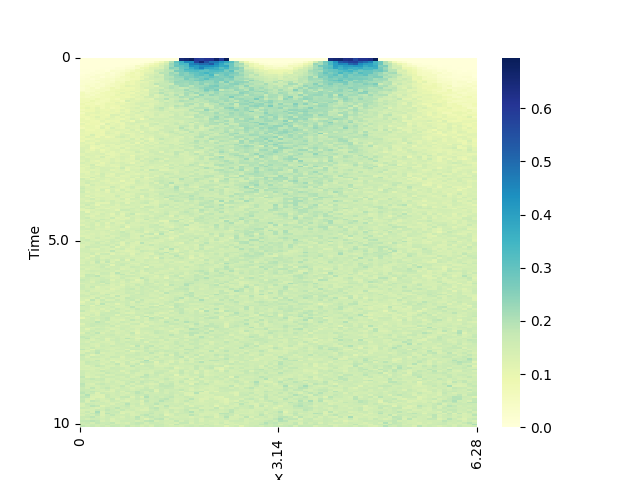}
        \caption{$\kappa=0.8$}
    \end{subfigure}
    \caption{Generalization on the randomness}
    \label{fig: Generalization}
\end{figure}

In this section, we show the generalization property of the neural network. For $\kappa=0.8, 2.5$ and fixed $N=3000$, we train the neural network with one generation of $\boldsymbol{\xi}_1$, and the results are illustrated in the previous sections. Then we run the algorithm with the trained neural network on a new generation $\boldsymbol{\xi}_2$, with all other parameters fixed. The results are illustrated in Figure~\ref{fig: Generalization}.

Since we consider feedback controls, our algorithm also generalizes well across different initial distributions. Specifically, training the neural network with one initial distribution enables it to provide accurate solutions even when tested with a different initial distribution.

\subsection{Systemic risk model}
\label{ssec:risk}
We consider a mean-field model introduced in \cite{carmona2013meanfieldgamessystemic} and computed by \cite{PW2}. Set $E=\R$, the representative agent follows the dynamics
$$
\d X_t=\left[\kappa\left(\mathbb{E}\left[X_t\right]-X_t\right)+\alpha_t\right] \mathrm{d} t+\sigma \d W_t, \quad X_0 \sim \mu_0,
$$
and the cost function is 
$$
J(\alpha)=\mathbb{E}\left[\int_0^T \tilde{f}\left(X_t, \mathbb{E}\left[X_t\right], \alpha_t\right) \mathrm{d} t+\tilde{g}\left(X_T, \mathbb{E}\left[X_T\right]\right)\right],
$$
where
$$
\tilde{f}(x, \bar{x}, a)=\frac{1}{2} a^2-q a(\bar{x}-x)+\frac{\eta}{2}(\bar{x}-x)^2, \quad \tilde{g}(x, \bar{x})=\frac{c}{2}(x-\bar{x})^2.
$$
In this example, for comparison we take the same parameters as in \cite{PW2}.  
Specifically, we take $\sigma=1$, $\kappa=0.6$, $q=0.8$, $T=0.2$, $C=2$, $\eta=2$, $\Delta t=0.2$, and 
solve the problem with $t=0$ and six different initial distributions considered in \cite{PW2}. The results are listed in the following table.

\begin{table}[h]
    \centering
    \begin{tabular}{|l|c|c|c|c|c|c|c|c|c|}
\hline \multicolumn{3}{|c|}{ Case 1 } & \multicolumn{3}{c|}{ Case 2 } & \multicolumn{3}{c|}{ Case 3 } \\
\hline Ours & \cite{PW2} & Analytic & Ours & \cite{PW2} & Analytic & Ours&  \cite{PW2} & Analytic \\
\hline 0.1649 & 0.1670 & 0.1642 & 0.1450 & 0.1495 & 0.1446 & 0.1449&0.1497 & 0.1446 \\

\hline  \multicolumn{3}{|c|}{ Case 4} & \multicolumn{3}{c|}{ Case 5} & \multicolumn{3}{c|}{ Case 6 } \\
\hline Ours & \cite{PW2} & Analytic & Ours & \cite{PW2} & Analytic & Ours&  \cite{PW2} & Analytic \\
\hline 0.1640 & 0.1675 & 0.1642 & 0.1810& 0.1824 & 0.1812 & 0.1775& 0.1792 & 0.1772 \\
\hline

\end{tabular}
    \caption{Numerical results of systematic risk model}
    \label{tab:my_label}
\end{table}

For our numerical algorithm, we employ a neural network with 2 hidden layers, each with 20 neurons, and take $m(K)=10$. The separating class consists of simple polynomials.  We train the 
neural network on a MacBook Pro and get the result in less than 20 minutes. 
The results \cite{PW2} were obtained on the graphic card GPU NVidia V100 32Go 
and in 3 days
(see Section 5, page 14 of \cite{PW2}).

{
\subsection{Crowd aversion planning problem}
In this example, we take $E:=\R^2$. Consider a  mean field control problem in $\R^2$ where
the representative particle chooses its drift, wishes to
avoid crowded regions, and is steered towards a prescribed \emph{target} law at
time $T$. Let $(X_t)_{t\in[0,T]}$ solve the controlled SDE
\begin{equation}\label{eq:crowd-dynamics}
  \d X_t = \alpha_t\,\d t + \sigma\,\d W_t, \qquad X_0 \sim \mu_0 \in \mathcal P_2(\R^2),
\end{equation}
where $(W_t)_{t\in [0,T]}$ is a Brownian motion, $\sigma>0$ is fixed, and
$\alpha\in \cA$ is an admissible control.
We denote by $\mu_t := \cL(X_t)$ the flow of marginal laws at time $t\in [0,T]$.

The running cost is defined by
\begin{equation}\label{eq:crowd-running-cost}
  \ell(x,\mu,a)
  = \frac12 |a|^2 + \lambda_c\,\Phi\big(K_h*\mu\big),
\end{equation}
where $\lambda_c>0$ is a crowd-aversion parameter, $K_h(x):=
\frac{1}{2\pi h^2}
\exp\!\Big(-\frac{|x|^2}{2h^2}\Big)$ is a Gaussian kernel, and $\Phi:\R_+\to\R_+$ is an increasing
function. In the numerical experiments we use
\begin{equation*}
  \Phi(r) = r^\gamma, \qquad \gamma\ge1,
\end{equation*}
so that the marginal cost of entering a crowded region increases 
with the local density.

At the terminal time the controller is penalized for the mismatch between the
terminal law $\mu_T$ and a prescribed target distribution
$\mu_{\mathrm{tar}}\in\cP_2(\R^2)$. We consider terminal cost 
$$
k\, G(\mu,\mu_{\mathrm{tar}}),
$$
where $k>0$ controls the strength of the terminal constraint, and a typical choice of $G$ is the quadratic Wasserstein-2 distance $\cW_2^2$.

Putting these ingredients together, the mean field control problem reads
\begin{equation}\label{eq:crowd-MFC-problem}
  V_k(\mu_0) := \inf_{\alpha\in \cA}
  \E\Bigg[
    \int_0^T \Big(
      \frac12 |\alpha_t|^2
      + \lambda_c\,\Phi\big(K_h*\mu_t(X_t)\big)
    \Big)\,
    \d t
    + k\,G\big(\mu_T,\mu_{\mathrm{tar}}\big)
  \Bigg].
\end{equation}
\begin{remark}
        As $h$ converges to zero, for every $x\in \R^2$, $K_h*\mu(x)$ converges to $\rho_{\mu}(x)$, where $\rho_{\mu}$ represents the density of $\mu$\,(if exists).
\end{remark}
In our numerical experiment, we take $\mu_{\mathrm{ini}}=\cN((0,0), 0.3^2I_2)$, $\mu_{\mathrm{tar}}:=\delta_{(1.5,1.5)}$, $\gamma=2.0$. Denote $x_{\mathrm{tar}}:=(1.5,1.5)$, then $\cW_2^2(\mu,\mu_{\mathrm{tar}})=\int_{\R^2}\|x-x_{\mathrm{tar}}\|^2\mu(\d x)$ for every measure $\mu\in \cP_2(\R^2)$.   We use a fully connected neural network with three hidden layers, each with $64$ neurons, and we feed the first eight moments of the empirical measure as measure features into the network. The network is trained with the Adam algorithm with learning rate $10^{-3}$, using $N=1024$ particles in the simulation for the training process.  For different pairs $(\lambda_c,k)$, we train the network separately and then visualize the resulting dynamics by simulating 30 particles under the learned optimal feedback control.
\begin{figure}
    \begin{subfigure}[b]
    {0.3\textwidth}
        \centering
        \includegraphics[width=\textwidth]{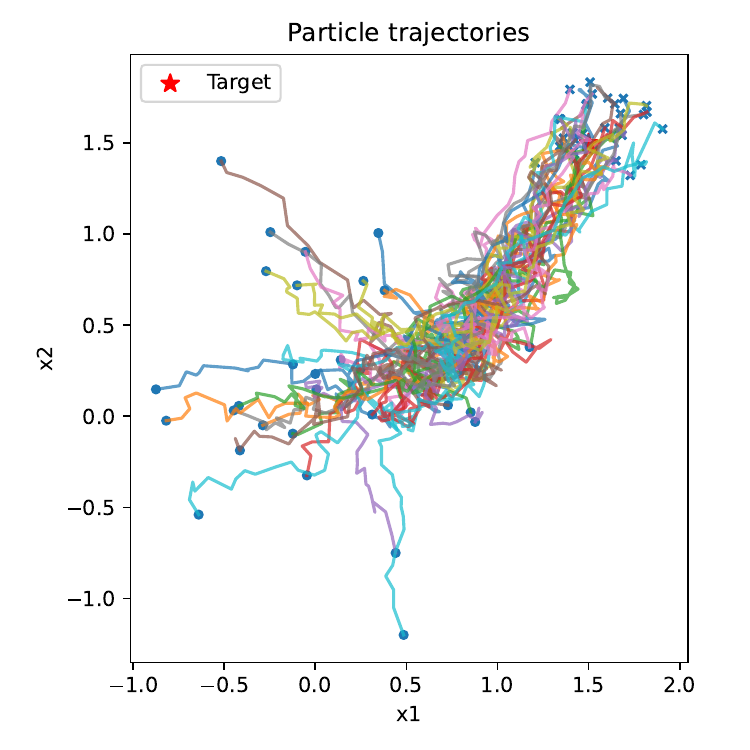}
        \caption{$\lambda_c=0,k=1000$}
        \label{fig: c0k1000}
    \end{subfigure}
    \begin{subfigure}[b]
    {0.3\textwidth}
        \centering
        \includegraphics[width=\textwidth]{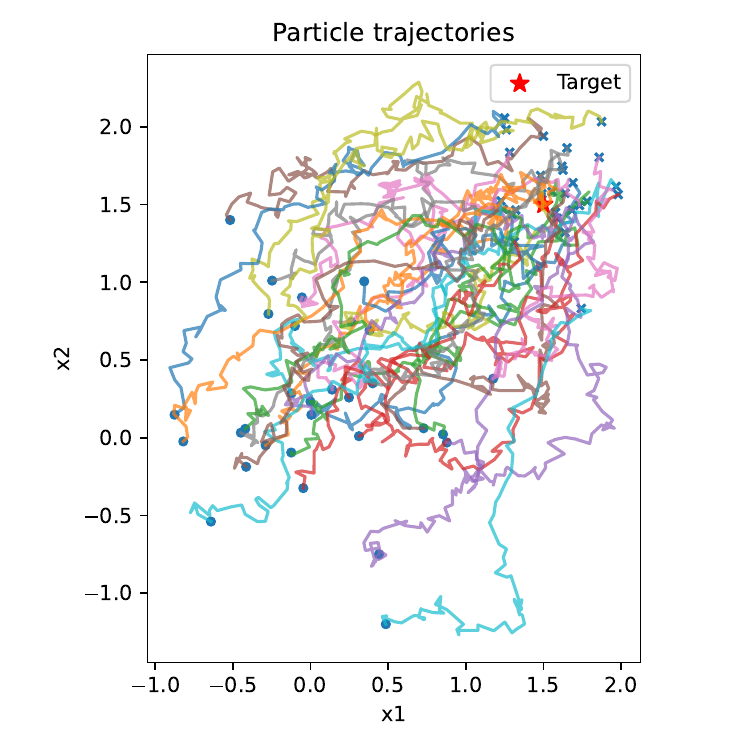}
        \caption{$\lambda_c=500,k=1000$}
        \label{fig:c500k1000}
    \end{subfigure}
    \begin{subfigure}[b]
    {0.3\textwidth}
        \centering
        \includegraphics[width=\textwidth]{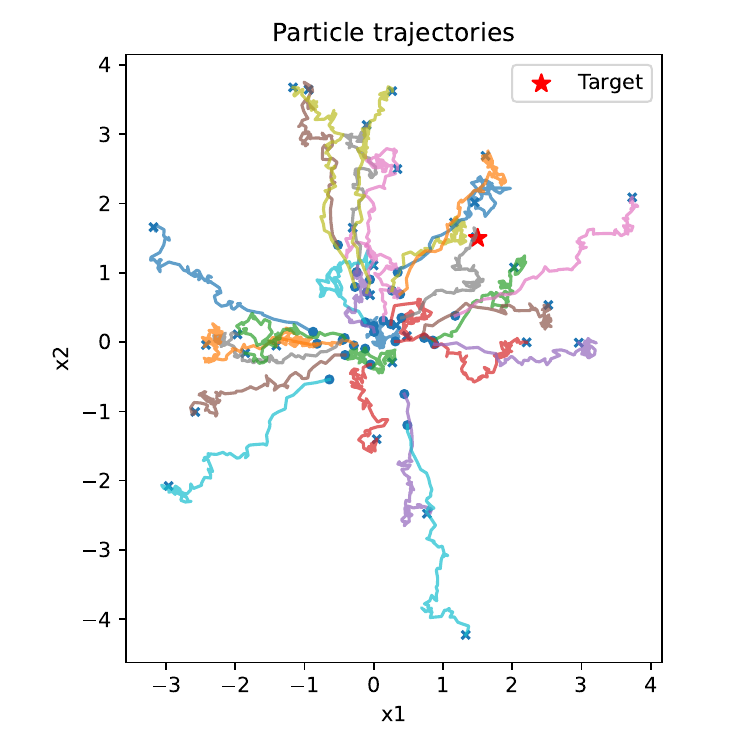}
        \caption{$\lambda_c=500,k=0$}
        \label{fig:c500k0}
    \end{subfigure}
     \caption{Numerical simulations of $30$ particles with different parameters.}
\end{figure}

From our numerical results, we can see from Figures~\ref{fig:c500k1000} and \ref{fig:c500k0}, for the same $\lambda_c$, increasing  $k$ will force the particles to move towards the target distribution, while Figures~\ref{fig: c0k1000} and \ref{fig:c500k1000} show that for fixed $k$, increasing $\lambda_c$ will encourage the particles to avoid high-density regions, leading to a more diffuse cloud and sparser optimal trajectories in the simulation.

\begin{remark}
    Beyond providing a simple testbed for our numerical algorithm, this
crowd-aversion control problem is also a stylized model for generative
modeling in high dimension when $\lambda_c=0$. This is closely related to soft-constraint Schrödinger bridge formulations of
generative AI \cite{garg2024softconstrainedschrodingerbridgestochastic,ma2025schrodingerbridgegenerativeai}. The efficiency of our algorithm shows that the same architecture and algorithm can in principle be applied
to high-dimensional feature spaces (e.g.\ embeddings of images or text),
providing a conceptually simple, control-theoretic counterpart to modern
diffusion and score-based generative models. We will leave the implementations of the generative model formulations for future research.
\end{remark}
}

\section{Convergence}
\label{sec:analysis}
In this section we prove that the value function obtained 
from our numerical algorithm converges to the mean field value function. Recall the definitions of 
value functions $v(\mu)$, $v^N(x^N)$, $v_K(\mu)$,
$v_K^N(x^N)$ and $V^N_K(X_0^N,\boldsymbol{\xi})$ from  
the previous sections.
Let $\mathcal{W}_1(\mu, \nu)$ be the Wasserstein one 
distance.

We impose the following regularity conditions on the given functions,
$$
b :\cX \times A \mapsto E,\ \
\sigma  :\cX \times A \mapsto \cM(d,m),\ \
\ell :\cX \times A \mapsto \R,\ \
G :E \times \cP_2(E)  \mapsto \R.
$$
\begin{assumption}[Standing Assumption]
\label{asm:main}
We assume that there are constants $c_*, c_K$  
and a function $c(\cdot)$ such that
for any $j \in \cT$, $K\in\N$, 
$(x,\mu), (y,\nu) \in E \times \cP_2(E)$, and $a,a'\in A$ with $|a-a'|\le 1$, 
\begin{itemize}
\item (Lipschitz Continuity) For $h =b,\sigma$,
\begin{align*}
    |h(j, x, \mu,a) - h(j,y, \nu,a')| &\leq c_* \left(|x - y| + \mathcal{W}_1(\mu, \nu) +|a-a'|\right),\\
    |G(x,\mu)-G(y,\nu)|&\leq c_* \left(|x - y| + \mathcal{W}_1(\mu, \nu)\right).
\end{align*}
\item (Local Lipschitz Continuity) 
$$\
|\ell(j, x, \mu,a) - \ell(j,y, \nu,a')| 
\leq c(|a|) \left(|x - y| + \mathcal{W}_1(\mu, \nu) +|a-a'|\right),
$$
\item (Growth Conditions) For  $h=b,\sigma$,
$|h(j,x,\mu,a)| \leq c(|a|) (1+|x|+\int_E |\tilde x|\mu(\d \tilde x))$, 
$$
0\le \ell(j,x,\mu,a) \le c_*(1+ |x|^2 + |a|^2+ \int_E |\tilde x|^2\mu(\d \tilde x)),
\quad 
0\le G(x,\mu) \le c_* (1+|x|^2+\int_E |\tilde x|^2\mu(\d \tilde x)).
$$
\item (Hypothesis class)
$$ 
|\alpha(j, x, \mu) - \alpha(j,y, \nu)|
\leq {{K}} \left(|x - y| + \mathcal{W}_1(\mu, \nu)\right),
\qquad 
|\alpha(i,x,\mu)| \le K,
\qquad \forall \alpha \in \cN_K.
$$

{\item(Parametric continuity) The map $\theta\mapsto \alpha(\cdot;\theta)$  is continuous from $\Theta_K$ into  $(\mathcal N_K,\|\cdot\|_\infty)$. }

\item (Noise Process)
The components of the noise
sequences  $\xi_j^{(i) }\in \R^m$ are independently sampled 
from a fixed measure $\gamma\in \cP(\R^m)$ 
with a compact support.
\end{itemize}
\end{assumption}

In particular, in all of our examples $\xi_j^{(i)}$ are time
discretizations  of a Brownian motion.  Additionally,
by combining the above Lipschitz estimates, we 
conclude that for each $K$ there is {{a constant 
$c_K$}} satisfying for $h=b,\sigma,\ell$,
\begin{align}
\label{eq:assume}
|h^\alpha(j, x, \mu) - h^\alpha(j,y, \nu)|
\leq c_K \left(|x - y| + \mathcal{W}_1(\mu, \nu)\right), \qquad \forall \alpha \in \cN_K. 
\end{align}
\begin{remark}
    Stronger Lipschitz continuity assumptions in $\cW_1$ are necessary for  \textit{almost sure} convergence in the propagation of chaos theorem~\ref{thm:poc}. Lipschitz continuity in $\cW_2$ only guarantees convergence in probability.
\end{remark}

\begin{lemma}[Stability with respect to the control]\label{lem:stab-alpha}
Under Assumption~\ref{asm:main}, there exists a constant $C_{K,T}>0$ such that
for all $\alpha,\beta\in\cN_K$ and all $N\in\N$,
$$
\sup_{j\in\cT}\cW_1(\mu_j^{N,\alpha},\mu_j^{N,\beta})
\le C_{K,T}\,\|\alpha-\beta\|_\infty,
\qquad
\sup_{j\in\cT}\cW_1(\mu_j^{\alpha},\mu_j^{\beta})
\le C_{K,T}\,\|\alpha-\beta\|_\infty,
$$
where $\|\alpha-\beta\|_\infty:=\sup_{j\in\cT}\sup_{x\in E,\mu\in \cP_2(E)}|\alpha(j,x,\mu)-\beta(j,x,\mu)|$.
\end{lemma}
\begin{proof}
    Synchronous coupling $X^\alpha,X^\beta$ with same noise gives
    \begin{align*}
        \mathbb E|X^\alpha_{j+1}-X^\beta_{j+1}|
&\le c_*\Big(\mathbb E|X^\alpha_j-X^\beta_j|+\mathcal W_1(\mu^\alpha_j,\mu^\beta_j)+\|\alpha-\beta\|_\infty\Big)\\
&\le c_*\Big(2\mathbb E|X^\alpha_j-X^\beta_j|+\|\alpha-\beta\|_\infty\Big).
    \end{align*}
  Discrete Grönwall concludes the proof. Same for the second inequality. 
\end{proof}

\begin{theorem}[Convergence of the algorithm]
\label{convergence}
    Suppose that the standing
    Assumption \ref{asm:main} and the universal approximation
    Assumption \ref{NN:UAT}  hold.  For a 
    given $\mu \in \cP_\infty(E)$ and any random initial data  
   $X^N_0=(X^{(1)}_0,...,X^{(N)}_0)$ 
   whose components are identically  and independently drawn from $\mu$,
    $$
    \lim_{K\rightarrow\infty}\lim_{N\rightarrow\infty} 
    |v(\mu)-V^N_K(X^N_0,\boldsymbol{\xi})|=0 \quad \text{a.s.}
    $$
\end{theorem}
{{
We complete the proof of the above theorem in the subsequent subsections.  
Below we outline the main steps. We fix $\eps>0$
and $\mu \in \cP_\infty(E)$.
\begin{itemize}
    \item In Corollary  \ref{Firstpart}, we show that
    there is $K_\eps$ such that
    $$
    |v(\mu)- v_{K_\eps}(\mu)| \le \eps.
    $$
    Recalling from \eqref{ee:vK} that $v_{K_\eps}(\mu)$
    is the minimal value over the hypothesis
    class $\cN_{K_\eps}$, we may restrict attention to feedback controls in $\cN_{K_\eps}$.
    In particular, elements of $\cN_{K_\eps}$ are
    Lipschitz continuous with a uniform bound $K_\eps$.
    \item  Using the uniform Lipschitz continuity of the controls in $\cN_{K_\eps}$, in Theorem~\ref{thm:poc} and Corollary~\ref{Cor:UPoC} we prove a uniform propagation of chaos result for 
    the corresponding payoffs.  
    \item The components of the
    particle system ${\boldsymbol{X}^{N,\alpha}}$ defined by \eqref{eq:N-dynamics}
    are not independent of each other. To overcome this difficulty,
    we leverage the previous step in subsection \ref{sec:radamacher}, 
    by defining an i.i.d.~approximation by replacing the 
    interaction  terms in \eqref{eq:N-dynamics} with the law
    of the process.  
    \item In Proposition \ref{prop:rclimit}, 
    we combine the Rademacher complexity results
    with the i.i.d.~approximation described above, 
    to prove that, for each $K>0$,
    $$
    \lim_{N\rightarrow\infty}  |v_K(\mu)-V^N_K(X^N_0,\boldsymbol{\xi})|=0\qquad \text{a.s.}.
    $$
\end{itemize} }}

\begin{remark}
\label{rem:asm}
To emphasize the novel steps in the above
convergence proof with clarity,
we chose to prove the above theorem under the  assumption that the supports of the noise process and the initial  condition are bounded. However, due to the global nature of 
the universal approximation Theorem \ref{thm:ourUAT},  tedious
but otherwise standard techniques extend this convergence
to noise and initial distributions that have polynomial decay.
\end{remark}

\begin{remark}
    The convergence in Theorem~\ref{convergence} is stated without an explicit rate. For the value functions of mean field control problems, explicit $N$-rates require
additional regularity of the coefficients and of the limiting value function.  Recent results in this direction include algebraic rates for McKean--Vlasov control values~\cite{cardaliaguet2023algebraic}, and sharp $1/N$-rates in regions of strong regularity~\cite{cardaliaguet2023sharpconvergenceratesmean}. We do not pursue such quantitative bounds here, since our goal is to prove convergence under assumptions compatible with the flexible architecture used in the numerical algorithm.
\end{remark}

\subsection{Neural Network approximation}

The following is essentially a direct consequence of the assumptions.

\begin{proposition}
\label{pr:close}
    For  every $\mu \in \cP_\infty(E)$,  $\eps>0$, and feedback  control $\alpha^* \in \cC$, 
    there exist  $K_\eps \in \N_+$ and $\alpha_\eps\in \cN_{K_\eps}$, 
    such that 
    $$
    |J(\mu,\alpha^*)-J(\mu,\alpha_\eps)|\leq \epsilon.
    $$
\end{proposition}

\begin{proof}
For any $\alpha \in \cC$, let $(X^\alpha_j)_{j\in \cT}$ denote 
the solutions of the  mean field dynamics \eqref{eq:mfdynamics} 
with the initial distribution $\mu$ and control 
$\alpha \in \cC$, and set
$\mu^\alpha_j:=\cL(X^\alpha_j)$.  
To simplify the notation, set $X_j:=X^{\alpha^*}_j$,
$\mu_j:=\mu^{\alpha^*}_j$.
 We complete the proof in several steps.
\vspace{5pt}

\noindent
{\emph{Step 1.}} (Set-up). Suppose that for $\alpha \in \cC$, 
there is a constant $c_0>0$ such that
$|\alpha(j,x,\eta)|  \le c_0$ for all $(j,x,\eta) \in \cX$.
As the support of $\mu$ is  bounded,
 in view of the growth assumption and the boundedness
 of the noise process $\xi$, there exists a non-decreasing function 
$r(c_0)$ depending only on $\mu$, 
 the maximum of $\xi$, and the function $c(\cdot)$ in the
 growth assumption such that
 \be
 \label{eq:est}
|\alpha(j,x,\eta)|  \le c_0, \quad
\forall (j,x,\eta) \in \cX,
\qquad \Rightarrow \qquad
 |X^\alpha_j| \le r(c_0), \qquad j \in \cT.
 \ee
In particular,
in view of  \eqref{eq:controlset}, $\alpha^*$ is bounded by $c(\alpha^*)$.
Hence, $ |X_j| \le r(c(\alpha^*))$ for all $j$. 
Define a compact subset $\mathfrakC$ of $\cX$ by
 $$
\mathfrakC:= \{ (j,x,\nu) \ :\ |x|\le  c_1, \ \
\text{support}(\nu) \subset B( c_1) \}, \quad
\text{where}\quad  c_1:= r(c(\alpha^*+1)),
$$
and $B(c_1)\subset E$ is a ball with radius $c_1$.
In view of Assumption \ref{NN:UAT}, for  $0<\delta\le 1$ to be chosen at the end of the proof,
there are $K_\delta \in \N_+$ and $\alpha_\delta\in \cN_{K_\delta}$, 
such that 
$$
|\alpha^*(j, x,\nu)-\alpha_\delta(j, x,\nu)|\leq \delta,
\qquad \forall\ (j,x,\nu) \in \mathfrakC.
$$
Then,  by \eqref{eq:controlset}, we have,
\be
\label{eq:bdd}
|\alpha_\delta(j,x,\nu)| \le c(\alpha^*)+1, 
\qquad \forall\ (j,x,\nu) \in \mathfrakC.
\ee

\noindent
{\emph{Step 2.}} (Support estimates).
Set $Y_j:= X^{\alpha_\delta}_j$, $\nu_j:= \mu^{\alpha_\delta}_j=\cL(Y_j)$.
We claim that $|Y_j| \le c_1$.  Indeed, define a feedback
control $\bar{\alpha} \in \cC$ by
$$
\bar{\alpha}(j,x,\eta):=  \alpha_\delta(j,x,\eta) \ 
\frac{\min\{|\alpha_\delta(j,x,\eta)|\ ; \ c(\alpha^*)+1\}}
{|\alpha_\delta(j,x,\eta)|}.
$$
Then, $|\bar{\alpha}(j,x,\eta)| \le c(\alpha^*)+1$ and by \eqref{eq:est} and 
the definition of $c_1$, we have
$$
|X^{\bar{\alpha}}_j| \le r(c(\alpha^*)+1) =c_1,
\qquad \Rightarrow \qquad
(j,X^{\bar{\alpha}}_j,\mu^{\bar{\alpha}}_j) \in \mathfrakC.
$$
By \eqref{eq:bdd},
$$
|\alpha_\delta(j,X^{\bar{\alpha}}_j,\mu^{\bar{\alpha}}_j)| \le c(\alpha^*)+1,
\qquad \Rightarrow \qquad
\alpha_\delta(j,X^{\bar{\alpha}}_j,\mu^{\bar{\alpha}}_j)
=\bar{\alpha}(j,X^{\bar{\alpha}}_j,\mu^{\bar{\alpha}}_j).
$$
The uniqueness of the mean field dynamics \eqref{eq:mfdynamics} 
imply that $Y_j= X^{\bar{\alpha}}_j$.  Hence, we conclude
that $|Y_j| =|X^{\bar{\alpha}}_j| \le c_1$.
In particular, we have shown that
$$
(j,X_j,\mu_j), (j, Y_j,\nu_j) \in \mathfrakC, \qquad \forall \ j \in \cT.
$$
Additionally, in view of the local Lipschitz assumption,
there is a constant $c_2$ \emph{independent of $\delta$}, 
such that all the relevant functions are Lipschitz with constant $c_2$.
\vspace{5pt}

\noindent
{\emph{Step 3.}} (Estimating $\E|X_j-Y_j|$).
In view of \eqref{eq:mfdynamics},
    \begin{align*}
        \E|X_{j+1}-Y_{j+1}|&\leq \E|b^{\alpha^*}(j,X_j,\mu_j)-b^{\alpha_\delta}(j,Y_j,\nu_j)|
        +\E|(\sigma^{\alpha^*}(j,X_j,\mu_j)-\sigma^{\alpha_\delta}(j,Y_j,\nu_j))\xi_j|\\
        &=: \cI +\cJ, \qquad j=0,\ldots, T-1.
    \end{align*}
    For the first term, the Lipschitz assumption on $b$ and
    \eqref{eq:controlset} implies that    
    \begin{align*}
        \cI&\leq \E|b^{\alpha^*}(j,X_j,\mu_j)-b^{\alpha^*}(j,Y_j,\nu_j)|
               +\E|(b^{\alpha^*}-b^{\alpha_\delta})(j,Y_j,\nu_j)|\\
        &\leq c_2 (\E|X_j-Y_j|+\cW_1(\mu_j,\nu_j)) + c_2 \E|(\alpha^*-\alpha_\delta)(j,Y_j,\nu_j)|.
    \end{align*}
    Since $\cW_1(\mu_j,\nu_j) \le \E|X_j-Y_j|$ and $(j,Y_j,\nu_j)\in \mathfrakC$,
    $\cI \le c_2 ( 2 \E|X_j-Y_j|+\delta)$.
    As the noise process is bounded, exactly the same argument applies to $\cJ$.  Hence,
    there is a constant $c_3>0$, again \emph{independent of $\delta$} such that
    $$
     \E|X_{j+1}-Y_{j+1}|\leq \cI +\cJ
     \le c_3 ( \E|X_j-Y_j|+\delta).
   $$     
     Starting from  $\E|X_0-Y_0|=0$, by induction we conclude that 
     there is a constant $c_4$ depending on $\mu,\alpha^*,T$ and
     the constants of Assumption \ref{asm:main}, but otherwise independent
     of $\delta$, satisfying
   $$
   \E|X_{j}-Y_{j}|\leq c_4\ \delta,\quad \forall\, j\in \cT.
   $$

\noindent
{\emph{Step 4.}} (Estimating the cost functional).
Since both $(j,X_j,\mu_j), (j,Y_j,\nu_j) \in \mathfrakC$,
the control processes  $\alpha^*(j,X_i,\mu_j)$
and $\alpha_\delta(j,Y_j,\nu_j)$ are bounded.
Then, the Lipschitz assumptions imply that
there is yet another constant $c_5$ \emph{independent
of $\delta$} such that 
$$
|J(\mu,\alpha^*) - J(\mu,\alpha_\delta) | \le c_5 \ \delta.
$$
We complete the proof by choosing $\delta = \epsilon/c_5$.
\end{proof} 

An immediate corollary is the following.
\begin{corollary}
\label{Firstpart}
    For every $\epsilon>0$ and $\mu \in \cP_\infty$, 
    there exists $K_\eps\in \mathbb N_+$, such that 
    $$
    |v(\mu)-v_{K_\eps}(\mu)|\leq \epsilon.
    $$
\end{corollary}
\begin{proof}
As $\cN_K\subset \cC$,  $v(\mu)\leq v_K(\mu)$. 
Choose $\alpha^*\in \cC$, such that
$J(\mu,\alpha^*)-\epsilon/2\leq v(\mu)$.
In view of Proposition \ref{pr:close}, there exists $K_\eps\in \N$ and $\alpha_\eps\in \cN_{K_\eps}$, such that
$
|J(\mu,\alpha^*)-J(\mu,\alpha_{\eps})|\leq \epsilon/2.
$
Hence,
$$
v_{K_\eps}(\mu)-\epsilon\leq J(\mu,\alpha_\eps)-\epsilon
\leq J(\mu,\alpha^*)-\epsilon/2\leq v(\mu)\leq v_{K_\eps}(\mu).
$$
\end{proof}

\subsection{Propagation of chaos}

The concept of \emph{propagation of chaos} provides a framework for understanding how, in the large-particle limit, systems of interacting particles simplify and become independent and identically distributed. 

Propagation of chaos was first introduced in \cite{Kac1956} in the context of kinetic theory, where the Boltzmann 
equation was derived as a limit of interacting particle systems. McKean extended this idea to systems governed 
by stochastic differential equations and established convergence to the McKean–Vlasov equation \cite{McKean1966}. 
Subsequent rigorous developments were provided by Funaki \cite{Funaki1984} and Sznitman \cite{Sznitman1991}, 
who provided a comprehensive treatment of propagation of chaos, including the convergence of empirical measures. 
Del Moral \cite{delmoral2004feynmankac} also discussed the discrete-time interaction system example 
and the propagation of chaos results in section 1.5.1. Similar results in discrete-time follow directly
from these studies.  {{To state our result, we
fix $K$ and $\cN_K$ and for each $j\in \cT$, define
$\Psi_j: \R^d \times \cP_2(\R^d) \times A \times \Xi \to \R^d$ by
$$
\Psi_j(x,\mu,a,\xi):=b(j,x,\mu,a)+\sigma(j,x,\mu,a)\,\xi.
$$
By Assumption~\ref{asm:main}, there exists a constant $\gamma > 0$,
depending on $K$, such that for all $j \in \cT$ and $\xi \in \Xi$:
\begin{equation*}
\big|\Psi_j(x,\mu,a,\xi) - \Psi_j(y,\nu,b,\xi)\big| 
\;\le\; 
\gamma \left( |x-y| + \cW_1(\mu, \nu) + |a-b| \right).
\end{equation*}
\begin{theorem}[Discrete-Time Propagation of Chaos]
\label{thm:poc}
For an arbitrary $\alpha \in \cN_K$,
let ${\boldsymbol{X}^{N,\alpha}}$ be the solution to the $N$-particle system \eqref{eq:N-dynamics},
$\mu^{N,\alpha}_j$ be its empirical measure,
and ${\boldsymbol{X}^\alpha}=(X_0,\ldots,X_T)$ be the solution of the McKean--Vlasov dynamics \eqref{eq:mfdynamics}.
Then,
$$
\lim_{N\rightarrow\infty} \sup_{j\in \cT} \cW_1\bigl(\mu^{N,\alpha}_j, \cL(X_j^\alpha)\bigr) = 0, 
\qquad \text{a.s.}
$$
\end{theorem}
\begin{proof}
{{We directly follow Section 8.3, 
    Theorem 8.3.2 and Theorem 8.3.3 in
    \cite{delmoral2004feynmankac}.}}
    
    Fix the control $\alpha \in \cN_K$ and the time horizon $T$. Recall $$\mu^{N,\alpha}_j=\frac{1}{N}\sum_{i=1}^N \
     \delta_{X^{(i),\alpha}_j},\quad \mu_j^{\alpha}=\cL(X^{\alpha}_j).$$
Since time is discrete and the horizon is finite, it is enough to prove that for each fixed $j \in \cT$,
$$
\lim_{N\to\infty} \cW_1\bigl(\mu_j^{N,\alpha}, \mu_j^\alpha\bigr) = 0 
\quad\text{a.s..}
$$

\medskip

\noindent\emph{Step 1. Coupling construction.}
We use a synchronous coupling construction. Let $(\xi_{j+1}^i)_{j\in\cT,\, i\ge 1}$ be an i.i.d.\ family of noise variables used in both systems.

For the mean field SDE,  let $(X_j^{(i),\alpha})_{j\in\cT}$, $1\le i\le N$, be $N$ i.i.d.\ copies of the McKean--Vlasov process solving
$$
X_{j+1}^{(i),\alpha} 
= \Psi_j\bigl(X_j^{(i),\alpha},\mu_j^\alpha,
\alpha(j,X_j^{(i),\alpha},\mu_j^\alpha),\xi_{j+1}^i\bigr),
$$
where $\mu_j^\alpha := \cL(X_j^{(i),\alpha})$ is the (deterministic) law of $X_j^\alpha$. 
Define the empirical measure of these i.i.d.\ limit particles by
$$
\bar{\mu}_j^{N,\alpha} := \frac{1}{N} \sum_{i=1}^N \delta_{X_j^{(i),\alpha}}.
$$
 
The interacting particles $X_j^{(i),N,\alpha}$ satisfy
$$
X_{j+1}^{(i),N,\alpha} 
= \Psi_j\bigl(X_j^{(i),N,\alpha},\mu_j^{N,\alpha},
\alpha(j,X_j^{(i),N,\alpha},\mu_j^{N,\alpha}),\xi_{j+1}^i\bigr),
$$
where $\mu_j^{N,\alpha} := \frac{1}{N} \sum_{i=1}^N \delta_{X_j^{(i),N,\alpha}}$ is the empirical measure of the interacting system. We assume identical initial conditions,
$$
X_0^{(i),N,\alpha} = X_0^{(i),\alpha}, \qquad i\ge1.
$$

All these processes are constructed on the same probability space, using the same noises.

\medskip

\noindent\emph{Step 2. Decomposition of the error.}
Fix $j\in\cT$. We decompose the discrepancy between $\mu_j^{N,\alpha}$ and $\mu_j^\alpha$ using an intermediate empirical measure:
\begin{equation}
\label{eq:triangle}
\cW_1(\mu_j^{N,\alpha}, \mu_j^\alpha) 
\;\le\; 
\underbrace{\cW_1(\mu_j^{N,\alpha}, \bar{\mu}_j^{N,\alpha})}_{\text{coupling error}} 
\;+\; 
\underbrace{\cW_1(\bar{\mu}_j^{N,\alpha}, \mu_j^\alpha)}_{\text{sampling error}}.
\end{equation}

For the sampling error,  the term $\cW_1(\bar{\mu}_j^{N,\alpha}, \mu_j^\alpha)$ is the Wasserstein distance between the empirical measure of i.i.d.\ variables and their common law. By Assumption~\ref{asm:main}, $\mu_j^\alpha$ has a finite first moment. Since the state space is Polish, by the strong law of large numbers for empirical measures in Wasserstein distance, we have
$$
\cW_1(\bar{\mu}_j^{N,\alpha}, \mu_j^\alpha) \longrightarrow 0,
\quad N\to\infty, \quad \text{a.s.}
$$
Let us denote
$$
\varepsilon_j^N := \cW_1(\bar{\mu}_j^{N,\alpha}, \mu_j^\alpha), 
$$
so that $\varepsilon_j^N \to 0$ a.s.\ as $N\to\infty$.

\medskip

\noindent\emph{Step 3. Propagation of the coupling error.}
We now control the first term in \eqref{eq:triangle}. For our specific coupling of $\mu_j^{N,\alpha}$ and $\bar{\mu}_j^{N,\alpha}$, the natural coupling pairs $X_j^{(i),N,\alpha}$ with $X_j^{(i),\alpha}$, $i=1,\dots,N$. Therefore
$$
\cW_1(\mu_j^{N,\alpha}, \bar{\mu}_j^{N,\alpha}) 
\;\le\; 
\frac{1}{N} \sum_{i=1}^N \big|X_j^{(i),N,\alpha} - X_j^{(i),\alpha}\big|.
$$
Define the average pathwise difference
$$
\Delta_j^N := \frac{1}{N} \sum_{i=1}^N \big|X_j^{(i),N,\alpha} - X_j^{(i),\alpha}\big| \;\ge 0.
$$
Then
$$
\cW_1(\mu_j^{N,\alpha}, \bar{\mu}_j^{N,\alpha}) \le \Delta_j^N.
$$

We now derive a recursive bound on $\Delta_j^N$. By Assumption \ref{asm:main}, the map $\Psi_j$ is Lipschitz in $(x,\mu,a)$ with constant $\gamma$ (with respect to $\cW_1$ in $\mu$), and $\alpha\in\cN_K$ is Lipschitz in $(x,\mu)$ with constant $c_K$. Thus, for each $i$ and each $j$,
\begin{align*}
&\big|X_{j+1}^{(i),N,\alpha} - X_{j+1}^{(i),\alpha}\big| \\
&= \Big|\Psi_j\bigl(X_j^{(i),N,\alpha},\mu_j^{N,\alpha},
\alpha(j,X_j^{(i),N,\alpha},\mu_j^{N,\alpha}),\xi_{j+1}^i\bigr)
-\Psi_j\bigl(X_j^{(i),\alpha},\mu_j^\alpha,
\alpha(j,X_j^{(i),\alpha},\mu_j^\alpha),\xi_{j+1}^i\bigr)\Big|
\\[0.2em]
&\le \gamma\Big(
\big|X_j^{(i),N,\alpha}-X_j^{(i),\alpha}\big|
+ \cW_1(\mu_j^{N,\alpha},\mu_j^\alpha)
+ \big|\alpha(j,X_j^{(i),N,\alpha},\mu_j^{N,\alpha})
      -\alpha(j,X_j^{(i),\alpha},\mu_j^\alpha)\big|
\Big)
\\
&\le \gamma\Big(
\big|X_j^{(i),N,\alpha}-X_j^{(i),\alpha}\big|
+ \cW_1(\mu_j^{N,\alpha},\mu_j^\alpha)
+ c_K\big(\big|X_j^{(i),N,\alpha}-X_j^{(i),\alpha}\big|
+ \cW_1(\mu_j^{N,\alpha},\mu_j^\alpha)\big)
\Big)
\\
&= \gamma(1+c_K)\Big(
\big|X_j^{(i),N,\alpha}-X_j^{(i),\alpha}\big|
+ \cW_1(\mu_j^{N,\alpha},\mu_j^\alpha)
\Big).
\end{align*}
Averaging over $i=1,\dots,N$ yields
$$
\Delta_{j+1}^N 
\le \gamma(1+c_K)\Big(\Delta_j^N + \cW_1(\mu_j^{N,\alpha},\mu_j^\alpha)\Big).
$$
Using \eqref{eq:triangle} and $\cW_1(\mu_j^{N,\alpha},\bar{\mu}_j^{N,\alpha})\le\Delta_j^N$, we have
$$
\cW_1(\mu_j^{N,\alpha},\mu_j^\alpha)
\le \Delta_j^N + \varepsilon_j^N.
$$
Substituting into the previous inequality gives the pathwise recursion
$$
\Delta_{j+1}^N 
\le \gamma(1+c_K)\Big(2\Delta_j^N + \varepsilon_j^N\Big).
$$

\medskip

\noindent\emph{Step 4. Inductive convergence in $j$.}
We now show by induction on $j$ that $\Delta_j^N \to 0$ a.s.\ as $N\to\infty$.

\emph{Base case $j=0$.}  
By construction, $X_0^{(i),N,\alpha} = X_0^{(i),\alpha}$ for all $i$, hence $\Delta_0^N=0$ for all $N$.

\emph{Inductive step.}  
Fix $j\in\{0,\dots,T-1\}$ and assume that $\Delta_j^N\to0$ a.s.\ as $N\to\infty$. We already know from Step~2 that $\varepsilon_j^N\to0$ a.s. Thus, almost surely, the real sequences $(\Delta_j^N)_N$ and $(\varepsilon_j^N)_N$ both converge to $0$, and the recursion
$$
\Delta_{j+1}^N 
\le \gamma(1+c_K)\Big(2\Delta_j^N + \varepsilon_j^N\Big),\quad \text{a.s.}
$$
implies
$$
\limsup_{N\to\infty} \Delta_{j+1}^N
\le \gamma(1+c_K)\big(2\cdot 0 + 0\big) = 0,\quad \text{a.s.}.
$$
Since $\Delta_{j+1}^N\ge0$, we conclude that $\Delta_{j+1}^N\to0$ as $N\to\infty$ almost surely.

Thus, by induction, for each $j\in\cT$ we have
$$
\Delta_j^N \longrightarrow 0,
\quad N\to\infty,\quad \text{a.s.}
$$

\medskip

\noindent\emph{Step 5. Conclusion in $\mathcal W_1$.}
Combining the decomposition \eqref{eq:triangle} with the bounds
$$
\cW_1(\mu_j^{N,\alpha},\bar{\mu}_j^{N,\alpha}) \le \Delta_j^N,
\qquad
\varepsilon_j^N = \cW_1(\bar{\mu}_j^{N,\alpha},\mu_j^\alpha),
$$
we obtain, for each fixed $j$,
$$
\cW_1(\mu_j^{N,\alpha},\mu_j^\alpha)
\le \Delta_j^N + \varepsilon_j^N
\longrightarrow 0,
\quad N\to\infty,\quad \text{a.s.}
$$
As discussed at the beginning of the proof, this implies
$$
\lim_{N\to\infty}\sup_{j\in\cT}
\cW_1(\mu_j^{N,\alpha},\mu_j^\alpha) = 0,
\quad \text{a.s.}
$$
\vspace{-10pt}
\end{proof}}

{\begin{corollary}[Uniform Propagation of Chaos]
\label{Cor:UPoC}
    Under Assumption~\ref{asm:main}, 
    $$
    \sup_{\alpha\in\mathcal N_K}\ \sup_{j\le T}\,\mathcal W_1(\mu_j^{N,\alpha},\mu_j^\alpha)\to 0 \qquad a.s.
    $$
    Consequently, for any bounded Lipschitz test function $\varphi$,
$$
    \lim_{N\rightarrow\infty}
   \ \sup_{\alpha \in \cN_K}\ 
    \sup_{j\in \cT}\left|\int_{\mathbb{R}^d} \varphi(x) 
    (\mu_j^{N,\alpha}(\d x) - \cL(X_j^\alpha)(\d x)) \right| =0,
    \quad \textit{a.s.}
$$
\end{corollary}
\begin{proof}
    Since $\Theta_K$ is bounded and finite-dimensional, for every $\eps>0$, there exists a finite $\varepsilon$-net  $\{\alpha^m\}_{m\le M}$:
     for every $\varepsilon>0$ there exist finitely many controls $\alpha^1,\dots,\alpha^M\in\mathcal N_K$  such that
$$
\forall \alpha\in\mathcal N_K\quad\exists\, m\le M:\quad \|\alpha-\alpha^m\|_\infty\le\varepsilon.$$
For any $\alpha$ pick $m$ with  $\|\alpha-\alpha^m\|_\infty\le\varepsilon$. Then,
\begin{align*}
    \sup_{j\in \cT}\mathcal W_1(\mu_j^{N,\alpha},\mu_j^\alpha)
&\le
\sup_{j\in \cT}\mathcal W_1(\mu_j^{N,\alpha},\mu_j^{N,\alpha^m})
+\sup_{j\in \cT}\mathcal W_1(\mu_j^{N,\alpha^m},\mu_j^{\alpha^m})
+\sup_{j\in \cT}\mathcal W_1(\mu_j^{\alpha^m},\mu_j^\alpha)\\
&\le
2C_{K,T}\varepsilon + \sup_{j\in \cT}\mathcal W_1(\mu_j^{N,\alpha^m},\mu_j^{\alpha^m}),
\end{align*}
where the second inequality comes from Lemma~\ref{lem:stab-alpha}. Take supremum over $\alpha\in \cN_K$ to get
$$
\sup_{\alpha\in \cN_K}\sup_{j\in \cT}\mathcal W_1(\mu_j^{N,\alpha},\mu_j^\alpha)
\le 2C_{K,T}\varepsilon+\max_{m\le M}\sup_{j\in \cT}\mathcal W_1(\mu_j^{N,\alpha^m},\mu_{j}^{\alpha^m}).
$$
By Theorem~\ref{thm:poc}, for each fixed $\alpha^m$, $\mathcal W_1(\mu_j^{N,\alpha^m},\mu_{j}^{\alpha^m})$ goes to $0$ a.s., 
and the maximum over finitely many $m$ goes to $0$ a.s.
We then let $\varepsilon\downarrow0$  to obtain
the first result. 

Let $\varphi$ be bounded and Lipschitz, with Lipschitz constant $L_\varphi$. By the Kantorovich--Rubinstein duality,
$$
\left|\int_{\R^d} \varphi(x)\,\big(\mu_j^{N,\alpha}-\mu_j^\alpha\big)(\d x)\right|
\le L_\varphi\,\cW_1(\mu_j^{N,\alpha},\mu_j^\alpha),
$$
for each $j,N$. Hence,
$$
\sup_{j\in\cT}
\left|\int_{\R^d} \varphi(x)\,\big(\mu_j^{N,\alpha}-\mu_j^\alpha\big)(\d x)\right|
\le L_\varphi\,\sup_{j\in\cT}\cW_1(\mu_j^{N,\alpha},\mu_j^\alpha),
$$
and the right-hand side converges to $0$ a.s.\ as $N\to\infty$ by the previous step.
\end{proof}}}

\subsection{Rademacher Complexity}
\label{sec:radamacher}

In this subsection, we follow the techniques used in \cite{RS} to bound the randomness and prove the uniform
law of large numbers. Consider initial data $\mu\in \cP_2(E)$ and  a random sequence 
$\boldsymbol{\xi}$.  As before we fix  a random initial data
$X_0=(X_0^{(1)},\ldots,X_0^{(N)})$ 
whose components are sampled i.i.d. from $\mu$.
Recall the notations $X_j(\bsx,X_0)$, $\mathfrakJ(X_0,\alpha,\bsx)$ in subsection \ref{sec:prob_formulation}, and 
$\mathfrak{J}^N(X_0,\alpha,\bsx)$
of \eqref{eq:jn}

For $i=1,\dots,N$ and $\alpha \in \cN_K$, let $\boldsymbol{Z}^{(i)}=X^\alpha(X^{(i)}_0,\bsx^{(i)})$ 
be the mean field dynamics from equation \eqref{eq:mfdynamics}
with the same control $\alpha\in \cN_K$.  As the noise processes 
$\bsx^{(i)}$ and $X_0^{(i)}$ are
all independent, $\boldsymbol{Z}=(\boldsymbol{Z}^{(1)},\dots, \boldsymbol{Z}^{(N)})$ 
are also identically and independently distributed, and $\cL(Z^{(i)}_j)= \cL(X_j)$,
where $X$ is the solution of  the mean field system \eqref{eq:mfdynamics}.  
Set
$$
w^N(X_0^N, \alpha,\boldsymbol{\xi}):=\frac{1}{N}\sum_{i=1}^N \left[\sum_{j=0}^{T-1} 
\ell^\alpha(j,Z^{(i)}_j,\cL(X_j))+G(\cL(X_T))\right]
=\frac{1}{N}\sum_{i=1}^N \mathfrakJ(X_0^{(i)},\alpha,\boldsymbol{\xi}^{(i)}).
$$
As each $Z^{(i)}_j$ is identically distributed to $X_j$, in view of \eqref{eq:jequal},
\be
\label{eq:wN}
\E[w^N(X_0^N, \alpha, \bsx)] = \E[\mathfrakJ(X_0,\alpha,\bsx)] = J(\mu,\alpha).
\ee

We note that $w^N(X_0^N, \alpha,\bsx)$ differs from $\mathfrakJ^N(X_0^N, \alpha, \bsx)$
slightly.  The difference is in the equations satisfied by $Z$ and the $X^N$ processes.
For the convenience of the reader we recall them:
\begin{align}
\label{eq:ZX}
Z^{(i)}_{j+1}&=Z_j^{(i)}+b^\alpha(j,Z_j^{(i)},\cL(X_j))
    +\sigma^\alpha(j,Z_j^{(i)},\cL(X_j)\ \xi_j^{(i)}), \\
    \nonumber
    X_{j+1}^{(i)}&=X_j^{(i)}+b^\alpha(j,X_j^{(i)},\mu^N_j))
    +\sigma^\alpha(j,X_j^{(i)},\mu^N_j))\ \xi_j^{(i)},\quad\quad
    j=0,...,T-1, \,i=1,\ldots,N, 
\end{align}
where $X$ is the solution of  the mean field system \eqref{eq:mfdynamics}.  
Notice that in the equation, for $X$ the interaction is through
the empirical measure $\mu^N_j=\mu^N(X^N_j)$.  Therefore,
the components of $X$ are coupled and are not necessarily independent
of each other.  On the other hand, $Z$ equation uses the
law of the process $\cL(Z^{(i)}_j)=\cL(X_j)$, which
is the same for each $i$, and therefore, as stated earlier the components
of $Z$ are independent.  In our analysis, we use the propagation
of chaos to argue that for each $\alpha$, the processes
$Z$ and $X$ are asymptotically close.  We then employ the 
Rademacher complexity estimates to  prove that in the large $N$ limit,
$w(X_0^N, \alpha,\boldsymbol{\xi})$ approximates $\E[\mathfrakJ(X_0^N, \alpha,\bsx)]= J(\alpha,\mu)$.

We first restate several classical definitions and results in Rademacher complexity, \cite{MRT}.
Set  
    $$
    \cG_K:=\{J(\cdot, \alpha, \cdot), \alpha\in \cN_K\}.
    $$
    The \emph{empirical Rademacher complexity} of $\cG_K$ on $(X_0,\boldsymbol{\xi})$ is defined as
    \begin{align*}
    R(\cG_K,N, X_0^N, \boldsymbol{\xi})
    &:=\E_{\sigma}[\sup_{g \in \cG_K} \frac{1}{N} 
    \sum_{i=1}^N \sigma_i g(X_0^{(i)},\boldsymbol{\xi}^{(i)})|(X_0^N, \bsx)],\\
   & = \E_{\sigma}[\sup_{\alpha \in \cN_K} \frac{1}{N} 
    \sum_{i=1}^N \sigma_i J(X_0^{(i)},\alpha,\bsx^{(i)})|(X_0^N, \bsx)],
    \end{align*}
    where $(\sigma_i)_{i=1,...,N}$ are i.i.d. random variables with values $\pm 1$. 
    The \emph{Rademacher complexity} of $\cG_K$ is
    $$
    r(\cG_K,N):=\E_{X_0^N,\boldsymbol{\xi}}[R(\cG_K,X_0^N, \boldsymbol{\xi})],
    $$
    where the expectation is taken over $\{(X_0^{(1)},\xi^{(1)}),\ldots,(X_0^{(N)},\xi^{(N)})\}$.
    Using \eqref{eq:wN}, we set
    $$
    \cE(\cN_K,N, X_0^N, \boldsymbol{\xi}):=\sup_{\alpha\in \cN_K}
     |w^N(\alpha,X_0^N, \boldsymbol{\xi}) - \E[w^N(\alpha,X_0^N,\boldsymbol{\xi})]|
     = \sup_{\alpha\in \cN_K}
     |w^N(\alpha,X_0^N, \boldsymbol{\xi}) - J(\mu,\alpha)| .
    $$
It is clear that for each $\alpha \in \cC$, 
$ |w^N(\alpha,X_0^N, \boldsymbol{\xi}) - \E w^N(\alpha,X_0^N,\boldsymbol{\xi})|$
converges to zero almost surely as $N$ tends to infinity.
Rademacher complexities allow us  to obtain uniform in $\alpha \in \cN_K$  estimates.
In particular, we utilize the following classical estimates \cite{MRT}.  A discussion of them is given in
Section 8 of \cite{RS} and we restate them for the convenience of the reader.
\begin{proposition}[\cite{MRT}]
\label{prop:rc}
For $\cG_K:=\cG(\cN_K)$, there exists a constant $\bar{c}_K$, depending on $K$, such that 
 for every  $\delta\in (0,1)$,  with the functions, 
    $$
     c(\cN_K,N, \delta):= 2r(\cG_K,N)+2 \bar{c}_K\sqrt{\frac{\ln{2/\delta}}{2N}},\quad
    C_e(\cN_K,N,X_0^N,\boldsymbol{\xi},\delta):=2R(\cG_K,N,X_0^N,\boldsymbol{\xi})
    +6\bar{c}_K\sqrt{\frac{\ln{2\delta}}{N}}.
    $$
the following estimates hold with probability $1-\delta$,
    $$
    \cE(\cN_K,N,X_0^N, \boldsymbol{\xi})\leq c(\cN_K,N,\delta)\leq 
    C_e(\cN_K,N,X_0^N, \boldsymbol{\xi},\delta).
    $$
    Additionally,
\be
\label{prop:rclimit}
\lim_{N\rightarrow\infty}r(\cG_K,N)=0.
\ee
\end{proposition}
\begin{remark}
    Instead of considering a fixed initial point and different independent noise
    processes as in \cite{RS}, we consider the joint product space of the initial position and path space, 
    then apply Rademacher complexity and related results for the product distribution. In other words, 
    by considering $\mu:=\mu_{x_0}$, where $x_0\in E$, we recover the results as  in \cite{RS}.
\end{remark}
Now we will use these results to bound the randomness in the value function.
\begin{proposition}\label{thirdpart}
    For every  $K\in \N$, every  initial data $X_0^N$
    whose components are i.i.d. sampled from  $\mu$, and every
    noise process $\bsx$,
    $$
   \lim_{N \to \infty}
    |v_K(\mu)-V^N_K(X_0^N,\boldsymbol{\xi})| =0, \quad \textit{a.s.}
    $$
\end{proposition}

\begin{proof}
    For  $\boldsymbol{\xi}$ and $\alpha\in \cN_K$, let $X, Z$ be as in \eqref{eq:ZX}
    and set $\mu_j^N:=\mu^N(X^N_j)$.  Definitions of the value functions and 
    \eqref{eq:jequal} imply that
    \begin{align*}
         |V^N_K(X_0^N,\boldsymbol{\xi}) &- v_K(\mu)| 
       \le \sup_{\alpha \in \cN_K} |\mathfrakJ^N(X_0^N,\alpha,\bsx)-J(\mu,\alpha)|\\
       &\leq
         \sup_{\alpha \in \cN_K} (|\mathfrakJ^N(X_0^N,\alpha,\bsx)-w^N(X_0^N,\alpha,\boldsymbol{\xi})|
        +|w^N(X_0^N,\alpha,\bsx)-\E w^N(X_0^N,\alpha,\bsx)|)\\
        &\le  \sup_{\alpha \in \cN_K}
        ( |\mathfrakJ^N(X_0^N,\alpha,\bsx)-w^N(X_0^N,\alpha,\boldsymbol{\xi})|)
        + \cE(\cN_K,N, X_0^N, \boldsymbol{\xi}).
    \end{align*}
    In view of the definitions of $Z$ and $w^N$,
    \begin{align*}
        |\mathfrakJ^N(X_0^N,\alpha,\bsx)&-w^N(X_0^N,\alpha,\boldsymbol{\xi})|\\
       & =|  \sum_{j=0}^{T-1}\frac{1}{N}\sum_{i=1}^N 
        (\ell^\alpha(j,X_j^{(i)},\mu^N_j)-\ell^\alpha(j,Z_j^{(i)},\cL(X_j)))
        +G(\mu^N_T)-G(\cL(X_T))  |\\
        &\le\sum_{j=0}^{T-1}\frac{1}{N} | \sum_{i=1}^N 
        \ell^\alpha(j,X_j^{(i)},\mu^N_j)-\ell^\alpha(j,Z_j^{(i)},\cL(X_j))|
        +c_* \cW_1(\mu^N_T,\cL(X_T)).
    \end{align*}
    For $j\in\{0,\ldots,T-1\}$,
    \begin{align*}
       & \frac{1}{N} | \sum_{i=1}^N
       \ell^\alpha(j,X_j^{(i)},\mu^N_j)-\ell^\alpha(j,Z_j^{(i)},\cL(X_j))|\\
      &\hspace{10pT}   \le 
      \frac{1}{N}| \sum_{i=1}^N
      \ell^\alpha(j,X_j^{(i)},\mu^N_j)-\ell^\alpha(j,X_j^{(i)},\cL(X_j))|
        +\frac{1}{N}\sum_{i=1}^N
        | \ell^\alpha(j,X_j^{(i)},\cL(X_j))-\ell^\alpha(j,Z_j^{(i)},\cL(X_j))|.
    \end{align*}
    Let $c_K$ be as in \eqref{eq:assume}. Then,
    $$
     \frac{1}{N}\sum_{i=1}^N
      |\ell^\alpha(j,X_j^{(i)},\mu^N_j)-\ell^\alpha(j,X_j^{(i)},\cL(X_j))|
      \le  c_K \cW_1(\mu^N_j,\cL(X_j)).
    $$
    Set $Z^N_j:=(Z^{(1)}_j,\ldots,Z^{(N)}_j)$.  Then,
     \begin{align*}
    \frac{1}{N} | &\sum_{i=1}^N
        \ell^\alpha(j,X_j^{(i)},\cL(X_j))-\ell^\alpha(j,Z_j^{(i)},\cL(X_j))| 
        = |\int \ell^\alpha(j,x,\cL(X_j))\ (\mu^N_j(\d x)-\mu^N_j(Z^N_j))(\d x)|\\
       & \le |\int \ell^\alpha(j,x,\cL(X_j))\ (\mu^N_j-\cL(X_j))(\d x)|
       +|\int \ell^\alpha(j,x,\cL(X_j))\ (\cL(X_j)-\mu^N_j(Z^N_j))(\d x)|\\
       & \le c_K (\cW_1(\mu^N_j,\cL(X_j)) + \cW_1(\cL(X_j),\mu^N_j(Z^N_j)).
       \end{align*}
    Combining the above inequalities and using
     Assumption~\ref{asm:main} and Corollary~\ref{Cor:UPoC}, 
     we conclude that 
     there is $c_1$ depending  on $K$ and $T$ such that
     \be
     \label{eq:1}
     \sup_{\alpha \in \cN_K} ( |\mathfrakJ^N(X_0^N,\alpha,\bsx)-w^N(X_0^N,\alpha,\boldsymbol{\xi})| )
     \rightarrow 0,\qquad \textit{a.s.}.
     \ee

    We estimate $\cE(\cN_K,N, X_0^N, \boldsymbol{\xi})$
    using the Rademacher complexity. 
    In Proposition \ref{prop:rc}, we take $\delta_N:=2e^{-2N\epsilon^2}$, 
    then with probability $1-\delta_N$,
    $$
    \sup_{\alpha\in \cN_K}\cE(\cN_K,N, X_0^N, \boldsymbol{\xi})
    \leq 2r(\cG_K,N)+2\bar{c}_K \epsilon.
    $$
    Since $\sum_N\delta_N< \infty$, by the Borel-Cantelli Lemma and \eqref{prop:rclimit},
    \begin{align}\label{eq:2}
         \limsup_{N\rightarrow\infty} \cE(\cN_K,N, X_0^N, \boldsymbol{\xi})
         \leq 2 \bar{c}_K \epsilon,\quad \textit{a.s.}
    \end{align}
We combine the above limit with \eqref{eq:1} to arrive at
 \begin{align*}
  \limsup_{N\rightarrow\infty} &|V_K^N(X_0^N,\boldsymbol{\xi})-v_K(\mu)|\\
  &\hspace{20pt}\le
 \limsup_{N\rightarrow\infty} \sup_{\alpha\in \cN_K}    
 (|\mathfrakJ^N(X_0^N,\alpha,\boldsymbol{\xi})-J(\mu,\alpha)| +
 \cE(\cN_K,N, X_0^N, \boldsymbol{\xi}))  \leq2 \bar{c}_K \epsilon,\quad \textit{a.s.}   
 \end{align*}
\end{proof}

\subsection{Proof of Theorem \ref{convergence}}
\begin{proof}
By the triangle inequality,
\begin{align*}
     |v(\mu)-V^N_K(X_0^N,\boldsymbol{\xi})|& \leq
     |v(\mu)-v_K(\mu)|+|v_K(\mu)-V^N_K(X_0^N,\boldsymbol{\xi})|.   
\end{align*}
Propositions \ref{Firstpart}  and  \ref{thirdpart} complete the proof.
\end{proof}

\section{Conclusions}
\label{sec:conclusion}
In this paper, we introduce and analyze a  machine learning--based numerical algorithm solving mean field optimal control problems.  
Unlike many existing strategies that rely on multiple simulations of $N$-particle systems, our framework employs a single large-population trajectory to learn a feedback control that depends on time, each agent's state, and the distribution.
This design allows for efficient computation and naturally incorporates the permutation invariance property of the particle system within the control.

From a theoretical standpoint, we establish that, under suitable regularity conditions, the large-population system converges to the mean field limit.
The key is to leverage classical propagation-of-chaos arguments and bound the randomness by Rademacher complexity, thereby theoretically verifying the correctness of our algorithm.

Beyond potential refinements of the theoretical arguments (for instance, refining convergence rates), there are also many interesting avenues for numerical implementations. As mentioned in Remark~\ref{rmk: seprarating class}, in a more complicated model, one can leverage the learning ability of neural networks by training them to learn a suitable separating class of functions. Another direction is to explore alternative neural network architectures or training strategies that could further improve computational scalability and approximation accuracy in high-dimensional state spaces. We believe these directions open new horizons for the machine learning--based algorithms of mean field control problems, enabling scalable and rigorous solutions for increasingly complex models.

\medskip

\appendix
\section{Proof of Proposition \ref{viscosity}}
In this appendix, we provide a convergence result under general assumptions on the initial data.
\begin{proposition}\label{viscosity}
    For every $K>0$, and a sequence of initial conditions $x^N\in E^N$
    satisfying $\mu^N(x^N)$ converging to $\mu$ in $\cW_2$,
    $v^N_K(x^N)$ converges to $v_K(\mu)$.
    Consequently, for any random initial data 
    $X^N_0=(X^{(1)}_0,\dots,X^{(N)}_0)$ whose 
    components are independently drawn from
    $\mu\in \cP_2(E)$ we have
    $$
    \lim_{N\rightarrow\infty}|v_K(\mu)-v^N_K(X^N_0)|=0 \quad \text{a.s.}
    $$
\end{proposition}
The proof of the first statement, given in the Appendix,
uses dynamic  programming principle and techniques from
the viscosity theory \cite{FS}. The second part of the theorem is an immediate
consequence of the law of large numbers.

In this section, we omit the dependence on $K$ and write $v^N$ for $v^N_K$ and
$v$ for $v_K$. 

\subsection{Dynamic programming and Hamilton-Jacobi-Bellman equation}
We fix $K$ and  define the mean field optimal control problem
starting from any time $t \in \cT$  as
\begin{equation}\label{mfdvaluefunc}
 v(t,\mu):=\inf_{\alpha\in \cN_K} J(t, \mu, \alpha)
 :=\inf_{\alpha\in \cN_K} \sum_{j=t}^{T-1} 
 \E\left[\ell^\alpha(j,\cL(X_j),X_j)+G(X_T,\cL(X_T)) \right],
\end{equation}
where as before $(X_j)_{j\in\cT}$ follows the same dynamics~\eqref{eq:mfdynamics} with control $\alpha$ and $\cL(X_t)=\mu$. As we consider feedback controls, by an abuse of notation, denote
$$
\ell^\alpha(j,\cL(X_j)):=\E[\ell^\alpha(j,\cL(X_j),X_j)],\qquad G(\cL(X_T)):=\E[G(X_T,\cL(X_T))].
$$
We have the following dynamic programming equation,
\begin{align}
\label{mf:dpp}
    v(t,\mu)&=\inf_{a\in \cN_K} \{\ell^\alpha(t,\mu)+v(t+1,\cL(X_{t+1}))\},\quad t={0,...,T-1},\\
    \nonumber
    v(T,\mu)&=G(\mu),
\end{align}
Similarly, we define the dynamic formulation of the $N$-particle system
\begin{align*}
    v^N(t,x^N)&:=\inf_{\alpha\in \cN_K}J^N(t,x^N,\alpha):=\inf_{\alpha\in \cN_K}\mathbb E\left[\sum_{j=t}^{T-1} \ell^{N,\alpha}(j,X^N _j)+G^N(X^N_T)\right],
    \quad \text{where}\\
    \ell^{N,\alpha}(j,x^N)&:=
    \frac{1}{N}
    \sum_{i=1}^N \ell^\alpha(j,x^{(i)},\mu^N(x^N)),
    \qquad \forall\ x^N=(x^{(1)},
    \ldots, x^{(N)}) \in E^N,
\end{align*}
and $(X^N_j)_{j\in \cT}$ follows the dynamics~\eqref{eq:N-dynamics} with control $\alpha$ and initial position $x^N\in E^N$.
Then by the dynamic programming principle,
\begin{align}\label{N-DPP}
    v^N(t,x^N)&=\inf_{\alpha\in \cN_K}\{\ell^{N,\alpha}(t,x^N)+\mathbb \E[v^N(t+1,X^N_{t+1})]\},\quad t={0,...,T-1},\\
    \nonumber
    v^N(T,x^N)&=G^N(x^N).
\end{align}
Given our standing assumptions~\ref{asm:main}, the following theorem is standard in mean field control and stochastic optimal control literature, similar results can be found in \cite{SY} and \cite{FS}.
\begin{theorem}\label{uniqvis}
Under the previous assumptions, the value function $v:\mathcal T\times \mathcal{P}_2(E)\rightarrow \mathbb R$ satisfies the dynamic programming principle~\eqref{mf:dpp}, and the value function $v^N:\cT \times E^N\rightarrow \mathbb R$ satisfies the dynamic programming principle~\eqref{N-DPP}.
\end{theorem}

\subsection{Barles-Perthame relaxed limits}
Following the classical technique \cite{BS},
for $\mu\in \mathcal{P}_2(E)$, we 
define the following relaxed limits by
\begin{align*}
    v^*(t,\mu)&:=\lim_{M \to \infty\ \delta\rightarrow 0}
    \sup\{ v^N(t,x^N):   N\ge M, \cW_2(\mu^N(x^N),\mu)\leq \delta\},\\
    v_{*}(t,\mu)&:=\lim_{M \to \infty\ \delta\rightarrow 0}
    \inf\{ v^N(t,x^N):   N\ge M,  \cW_2(\mu^N(x^N),\mu)\leq \delta\}.
\end{align*}
Clearly, $v^*(\mu)$ is upper semi-continuous, 
and $v_{*}(\mu)$ is lower semi-continuous in the Wasserstein two metric.
  
\begin{proposition}\label{prop:subsolu}
    For $t\in \cT,\mu\in \cP(E)$, $v^*(T,\mu)=v_*(T,\mu)=G(\mu)$ and for $t=0,\ldots,T-1$,
    \begin{align*}
        v^*(t,\mu)&\leq \inf_{\alpha\in \cN_K} \{\ell^\alpha(t,\mu)+v^*(t+1,\cL(X_{t+1}))\},\\
        v_*(t,\mu)&\geq \inf_{\alpha\in \cN_K} \{\ell^\alpha(t,\mu)+v_*(t+1,\cL(X_{t+1}))\}.
    \end{align*}
\end{proposition}
\begin{proof}
    By the definitions of $v^N$, 
    $v^N(T,x^N)=G(\mu^N(x^N))=v(T,\mu^N(x^N))$. 
    Since, by Assumption \ref{asm:main},
    $G$ is Lipschitz continuous and by $\cW_2(\mu^N(x^N),\mu)\rightarrow 0$,
    $$
    v^*(T,\mu)=v_*(T,\mu)=G(\mu).
    $$
    We now use \eqref{N-DPP} and the continuity of $\ell$
    to arrive at
    \begin{align*}
    v^*(t,\mu)&= \lim_{M \to \infty\ \delta\rightarrow 0}
    \sup \{ v^N(t,x^N):  N\ge M,  \cW_2(\mu^N(x^N),\mu)\leq \delta\}\\
    & =\lim_{M \to \infty\ \delta\rightarrow 0}
     \sup \{  \inf_{\alpha\in \cN_K} \{\ell^{N,\alpha}(t,x^N)+\E v^N(t+1,X^N_{t+1})\}\ :\ 
     N\ge M, \cW_2(\mu^N(x^N),\mu)\leq \delta\})\\
     & \le \inf_{\alpha\in \cN_K} 
     \{ \lim_{M \to \infty\ \delta\rightarrow 0}
     \sup\{ \ell^{N,\alpha}(t,x^N)+\E v^N(t+1,X^N_{t+1})\}\ :\ 
    N\ge M,  \cW_2(\mu^N(x^N),\mu)\leq \delta\})\\
       & = \inf_{\alpha\in \cN_K}
     \{ \ell^{\alpha}(t,\mu)+\E v^*(t+1,\cL(X_{t+1}))\},
    \end{align*}
    where in the final equality we used the facts that,
    since $\mu^N(x^N)$ tends to $\mu$,  $\mu^N(X_{t+1})$ converges to $\cL(X_{t+1})$,
    and $v^N(t+1,X^N_{t+1})$ is a function of $\mu^N(X_{t+1})$.  
   
    A similar argument shows that $v_*$ is a supersolution of \ref{mf:dpp}.
\end{proof}
\subsection{Proof of Proposition \ref{viscosity}}
\begin{proof}
    We prove by backward induction. By Proposition \ref{prop:subsolu}, 
    $v_*(T,\mu)=v^*(T,\mu)=v(T,\mu)=G(\mu)$.
    Suppose that $v^*(t,\mu)=v^*(t,\mu)=v(t,\mu)$ for some  $t\in \cT$.
    Then, by the dynamic programming inequalities,
    \begin{align*}
    v^*(t-1,\mu) &\leq \inf_{\alpha\in \cN_K} \{\ell^\alpha(t-1,\mu) + v^*(t,\cL(X_t))\}
    = \inf_{\alpha\in \cN_K} \{\ell^\alpha(t-1,\mu)+ v(t,\cL(X_t))\} = v(t-1,\mu)\\
    &= \inf_{\alpha\in \cN_K} \{\ell^\alpha(t-1,\mu)+ v_*(t,\cL(X_t))\} \le v_*(t-1,\mu).
    \end{align*}
    Since  by definition, we always have
    $v_*(t-1,\mu)\leq v(t-1,\mu) \le v^*(t-1,\mu)$.  We conclude that
  $v_*(t-1,\mu)= v(t-1,\mu) = v^*(t-1,\mu)$.
  Thus by induction,
      $$
    \lim_{N\rightarrow\infty}|v(\mu)-v^N(x^N)|=0.
    $$
\end{proof}

\bibliography{MFC-Numerics}

@book{grohs2023proof,
  title={A proof that artificial neural networks overcome the curse of dimensionality in the numerical approximation of Black--Scholes partial differential equations},
  author={Grohs, Philipp and Hornung, Fabian and Jentzen, Arnulf and Von Wurstemberger, Philippe},
  volume={284},
  number={1410},
  year={2023},
  publisher={American Mathematical Society}
}

@article{becker2019deep,
  title={Deep optimal stopping},
  author={Becker, Sebastian and Cheridito, Patrick and Jentzen, Arnulf},
  journal={Journal of Machine Learning Research},
  volume={20},
  number={74},
  pages={1--25},
  year={2019}
}

@article{beck2020overview,
  title={An overview on deep learning-based approximation methods for partial differential equations},
  author={Beck, Christian and Hutzenthaler, Martin and Jentzen, Arnulf and Kuckuck, Benno},
 journal = {Discrete and Continuous Dynamical Systems - B},
volume = {28},
number = {6},
pages = {3697-3746},
year = {2023}
}

@book{MRT,
  title={Foundations of machine learning},
  author={Mohri, M. and Rostamizadeh, A. and Talwalkar, A.},
  year={2018},
  publisher={MIT press}
}

@article{BS,
  title={Convergence of approximation schemes for fully nonlinear second order equations},
  author={Barles, Guy and Souganidis, Panagiotis E},
  journal={Asymptotic analysis},
  volume={4},
  number={3},
  pages={271--283},
  year={1991}
}

@article{BP,
author={Barles,G. and Perthame,B.},
year={1988},
title={Exit Time Problems in Optimal Control and Vanishing Viscosity Method},
journal={SIAM Journal on Control and Optimization},
volume={26},
number={5},
pages={1133-16},
}

@article{HanE,
author = {Han, Jiequn and E, Weinan},
year = {2016},
month = {11},
pages = {},
title = {Deep Learning Approximation for Stochastic Control Problems}, 
notes = {In \textit{Deep Reinforcement Learning Workshop}}, 
journal = {NIPS}
}

@article {HanJentzenE,
	author = {Han, Jiequn and Jentzen, Arnulf and E, Weinan},
	title = {Solving high-dimensional partial differential equations using deep learning},
	volume = {115},
	number = {34},
	pages = {8505--8510},
	year = {2018},
	publisher = {National Academy of Sciences},
	journal = {Proceedings of the National Academy of Sciences}
}

@article{RS,
      title={Deep Empirical Risk Minimization in finance: looking into the future}, 
      author={A. Max Reppen and H. Mete Soner},
      year={2023},
      volume={33},
  number={1},
   journal={ Mathematical Finance},
  pages={116--145},
}

@article{RST,
      title={Deep Stochastic Optimization in Finance}, 
      author={A. Max Reppen and H. Mete Soner and Valentin Tissot-Daguette},
     journal={Digital Finance},
  volume={5},
  number={1},
  pages={91--111},
  year={2023}
}

@inproceedings{huang_nash_2007,
	title = {The {Nash} certainty equivalence principle and {McKean}-{Vlasov} systems: {An} invariance principle and entry adaptation},
	booktitle = {46th {IEEE} {Conference} on {Decision} and {Control}},
	author = {Huang, Minyi and Caines, Peter E. and Malham\'{e}, Roland P.},
	year = {2007},
	pages = {121--126}
}

@article{huang_large_2006,
	title = {Large population stochastic dynamic games: closed-loop {McKean}-{Vlasov} systems and the {Nash} certainty equivalence principle},
	volume = {6},
	number = {3},
	journal = {Communications in Information \& Systems},
	author = {Huang, Minyi and Malham\'{e}, Roland P. and Caines, Peter E.},
	year = {2006},
	pages = {221--252}
}

@inproceedings{huang_individual_2003,
	title = {Individual and mass behaviour in large population stochastic wireless power control problems: centralized and {Nash} equilibrium solutions},
	volume = {1},
	booktitle = {42nd {IEEE} {International} {Conference} on {Decision} and {Control}},
	author = {Huang, Minyi and PE, Caines and Malham\'{e}, R.P.},
	year = {2003},
	pages = {98--103}
}

@article{lasry_mean_2007,
	title = {Mean field games},
	volume = {2},
	number = {1},
	journal = {Japanese Journal of Mathematics},
	author = {Lasry, Jean-Michel and Lions, Pierre-Louis},
	year = {2007},	
	pages = {229--260}
}

@article{lasry_jeux_2006,
	title = {Jeux \`{a} champ moyen. {II} – {Horizon} fini et contr\^{o}le optimal},
	volume = {343},
	number = {10},
	journal = {Comptes Rendus Mathematique},
	author = {Lasry, Jean-Michel and Lions, Pierre-Louis},
	year = {2006},
	pages = {679--684}
}

@article{lasry_jeux_2006-1,
	title = {Jeux \`{a} champ moyen. {I} – {Le} cas stationnaire},
	volume = {343},
	number = {9},
	journal = {Comptes Rendus Mathematique},
	author = {Lasry, Jean-Michel and Lions, Pierre-Louis},
	year = {2006},
	pages = {619--625}
}

@article{CL1,
  title={Convergence analysis of machine learning algorithms for the numerical solution of mean field control and games {I}: The ergodic case},
  author={Carmona, Ren{\'e} and Lauri{\`e}re, Mathieu},
  journal={SIAM Journal on Numerical Analysis},
  volume={59},
  number={3},
  pages={1455--1485},
  year={2021},
  publisher={SIAM}
}

@article{CL2,
  title={Convergence analysis of machine learning algorithms for the numerical solution of mean field control and games: {II}—the finite horizon case},
  author={Carmona, Ren{\'e} and Lauri{\`e}re, Mathieu},
  journal={The Annals of Applied Probability},
  volume={32},
  number={6},
  pages={4065--4105},
  year={2022},
  publisher={Institute of Mathematical Statistics}
}

@inbook{CL3, 
author={Carmona, Ren{\'e} and Lauri{\`e}re, Mathieu}, 
title={Deep Learning for Mean Field Games and Mean Field Control with Applications to Finance}, 
booktitle={Machine Learning and Data Sciences for Financial Markets: A Guide to Contemporary Practices}, 
publisher={Cambridge University Press}, year={2023}, 
pages={369–392}}

@book{FS,
  title={Controlled Markov processes and viscosity solutions},
  author={Fleming, Wendell and Soner, H.\ Mete},
  volume={25},
  year={2006},
  publisher={Springer Science \& Business Media}
}

@book{CD,
  title={Probabilistic theory of mean field games with applications I-II},
  author={Carmona, Ren{\'e} and Delarue, Fran{\c{c}}ois and others},
  year={2018},
  publisher={Springer}
}

@article{Osher,
  title={A machine learning framework for solving high-dimensional mean field game and mean field control problems},
  author={Ruthotto, Lars and Osher, Stanley J and Li, Wuchen and Nurbekyan, Levon and Fung, Samy Wu},
  journal={Proceedings of the National Academy of Sciences},
  volume={117},
  number={17},
  pages={9183--9193},
  year={2020},
  publisher={National Academy of Sciences}
}

@article{reisinger,
  title={A fast iterative PDE-based algorithm for feedback controls of nonsmooth mean-field control problems},
  author={Reisinger, Christoph and Stockinger, Wolfgang and Zhang, Yufei},
  journal={SIAM Journal on Scientific Computing},
  volume={46},
  number={4},
  pages={A2737--A2773},
  year={2024},
  publisher={SIAM}
}

@article{HHL,
  title={Learning high-dimensional {M}c{K}ean--{V}lasov forward-backward stochastic differential equations with general distribution dependence},
  author={Han, Jiequn and Hu, Ruimeng and Long, Jihao},
  journal={SIAM Journal on Numerical Analysis},
  volume={62},
  number={1},
  pages={1--24},
  year={2024},
  publisher={SIAM}
}

@article{PhamConv,
author = {Bachouch, Achref and Huré, Côme and Pham, Huyên and Langrené, Nicolas},
title = {Deep Neural Networks Algorithms for Stochastic Control Problems on Finite Horizon: {C}onvergence Analysis},
journal = {SIAM Journal on Numerical Analysis},
volume = {59},
number = {1},
pages = {525-557},
year = {2021}
}

@Article{PhamNum,
  author={Achref Bachouch and Côme Huré and Nicolas Langrené and Huyên Pham},
  title={{Deep Neural Networks Algorithms for Stochastic Control Problems on Finite Horizon: {N}umerical Applications}},
  journal={Methodology and Computing in Applied Probability},
  year=2022,
  volume={24},
  number={1},
  pages={143-178},
  month={March}
}

@article{PhamActorCritic,
  title={Actor-critic learning for mean-field control in continuous time},
  author={Frikha, Noufel and Germain, Maximilien and Lauri{\`e}re, Mathieu and Pham, Huy{\^e}n and Song, Xuanye},
  journal={arXiv:2303.06993},
  year={2023}
}

@article{PW1,
  title={Actor-critic learning algorithms for mean-field control with moment neural networks},
  author={Pham, Huy{\^e}n and Warin, Xavier},
  journal={Methodology and Computing in Applied Probability},
  volume={27},
   number={1},
  pages={13},
  year={2025}  
}

@article{PW2,
  title={Mean-field neural networks-based algorithms for {M}c{K}ean-{V}lasov control problems},
  author={Pham, Huy{\^e}n and Warin, Xavier},
  journal={arXiv:2212.11518},
  year={2022}
}

@article{PW_wasserstein,
author = {Pham, Huyên and Warin, Xavier},
title = {Mean-field neural networks: Learning mappings on {W}asserstein space},
journal={Neural Networks},
  volume={168},
  pages={380--393},
  year={2023},
  publisher={Elsevier}
}

@article{CCS,
  title={Synchronization in a {K}uramoto mean field game},
  author={Carmona, Rene and Cormier, Quentin and Soner, H Mete},
  journal={Communications in Partial Differential Equations},
  volume={48},
  number={9},
  pages={1214--1244},
  year={2023}
}

@misc{HS,
      title={Potential Mean-Field Games and Gradient Flows}, 
      author={Felix Höfer and H. Mete Soner},
     journal={Stochastic Processes and their Applications},
  pages={104971},
  year={2026}
      }

@inproceedings{Kac1956,
  author = {Mark Kac},
  title = {Foundations of Kinetic Theory},
  booktitle = {Proceedings of the Third Berkeley Symposium on Mathematical Statistics and Probability},
  volume = {3},
  pages = {171--197},
  year = {1956}
}

@article{McKean1966,
  author = {Henry P. McKean},
  title = {A Class of Markov Processes Associated with Nonlinear Parabolic Equations},
  journal = {Proceedings of the National Academy of Sciences},
  volume = {56},
  number = {6},
  pages = {1907--1911},
  year = {1966}
}

@article{Funaki1984,
  author = {Toshio Funaki},
  title = {A Certain Class of Diffusion Processes Associated with Nonlinear Parabolic Equations},
  journal = {Zeitschrift für Wahrscheinlichkeitstheorie und Verwandte Gebiete},
  volume = {67},
  pages = {331--348},
  year = {1984}
}

@article{lacker,
  title={A general characterization of the mean field limit for stochastic differential games},
  author={Lacker, Daniel},
  journal={Probability Theory and Related Fields},
  volume={165},
  pages={581--648},
  year={2016},
  publisher={Springer}
}

@article{SY,
  title={Viscosity Solutions of the {E}ikonal Equation on the 
{W}asserstein space},
  author={Soner, H Mete and Yan, Qinxin},
  journal={Applied Mathematics \& Optimization},
  volume={90},
  number={1},
  pages={1},
  year={2024},
  publisher={Springer}
}

@article{SY1,
  title={Viscosity solutions for {M}c{K}ean--{V}lasov control on a torus},
  author={Soner, H Mete and Yan, Qinxin},
  journal={SIAM Journal on Control and Optimization},
  volume={62},
  number={2},
  pages={903--923},
  year={2024}
}

@book{Sznitman1991,
  author = {Alain-Sol Sznitman},
  title = {Topics in Propagation of Chaos},
  booktitle = {École d’Été de Probabilités de Saint-Flour XIX—1989},
  publisher = {Springer, Berlin, Heidelberg},
  pages = {165--251},
  year = {1991}
}

@misc{carmona2013meanfieldgamessystemic,
      title={Mean Field Games and Systemic Risk}, 
      author={Rene Carmona and Jean-Pierre Fouque and Li-Hsien Sun},
      year={2013},
      eprint={1308.2172},
      archivePrefix={arXiv},
      primaryClass={q-fin.PR},
      url={https://arxiv.org/abs/1308.2172}, 
}

@article{DH,
  title={Deep hedging},
  author={Buehler, Hans and Gonon, Lukas and Teichmann, Josef and Wood, Ben},
  journal={Quantitative Finance},
  volume={19},
  number={8},
  pages={1271--1291},
  year={2019}
}

@book{delmoral2004feynmankac,
    title = {Feynman-Kac Formulae: Genealogical and Interacting Particle Systems with Applications},
    author = {Del Moral, Pierre},
    year = {2004},
    publisher = {Springer Science \& Business Media}
}

@misc{CST:2025,
      title={Global universal approximation of functional input maps on weighted spaces}, 
      author={Christa Cuchiero and Philipp Schmocker and Josef Teichmann},
      year={2025},
      eprint={2306.03303},
      archivePrefix={arXiv},
      primaryClass={stat.ML},
      url={https://arxiv.org/abs/2306.03303}, 
}

@article{DLZ,
  title={Deep learning for population-dependent controls in mean field control problems with common noise},
  author={Dayanikli, Gokce and Lauri{\`e}re, Mathieu and Zhang, Jiacheng},
  journal={arXiv:2306.04788},
  year={2023}
}

@misc{garg2024softconstrainedschrodingerbridgestochastic,
      title={Soft-constrained Schrodinger Bridge: a Stochastic Control Approach}, 
      author={Jhanvi Garg and Xianyang Zhang and Quan Zhou},
      year={2024},
      eprint={2403.01717},
      archivePrefix={arXiv},
      primaryClass={stat.ML},
      url={https://arxiv.org/abs/2403.01717}, 
}

@misc{ma2025schrodingerbridgegenerativeai,
      title={Schr\"odinger bridge for generative AI: Soft-constrained formulation and convergence analysis}, 
      author={Jin Ma and Ying Tan and Renyuan Xu},
      year={2025},
      eprint={2510.11829},
      archivePrefix={arXiv},
      primaryClass={cs.LG},
      url={https://arxiv.org/abs/2510.11829}, 
}

@article{GLPW,
  title={DeepSets and their derivative networks for solving symmetric PDEs},
  author={Germain, Maximilien and Lauri{\`e}re, Mathieu and Pham, Huy{\^e}n and Warin, Xavier},
  journal={Journal of Scientific Computing},
  volume={91},
  number={2},
  pages={63},
  year={2022},
  publisher={Springer}
}

@article{Bahlali_2017,
   title={Existence and optimality conditions for relaxed mean-field stochastic control problems},
   volume={102},
   ISSN={0167-6911},
   url={http://dx.doi.org/10.1016/j.sysconle.2016.12.009},
   DOI={10.1016/j.sysconle.2016.12.009},
   journal={Systems \& Control Letters},
   publisher={Elsevier BV},
   author={Bahlali, Khaled and Mezerdi, Meriem and Mezerdi, Brahim},
   year={2017},
   month=Apr, pages={1–8} }

@article{yarotsky2017error,
  title={Error bounds for approximations with deep {ReLU} networks},
  author={Yarotsky, Dmitry},
  journal={Neural Networks},
  volume={94},
  pages={103--114},
  year={2017},
  publisher={Elsevier}
}

@article{petersen2018optimal,
  title={Optimal approximation of piecewise smooth functions using deep {ReLU} neural networks},
  author={Petersen, Philipp and Voigtlaender, Felix},
  journal={Neural Networks},
  volume={108},
  pages={296--330},
  year={2018},
  publisher={Elsevier}
}

@misc{cardaliaguet2023sharpconvergenceratesmean,
      title={Sharp convergence rates for mean field control in the region of strong regularity}, 
      author={Pierre Cardaliaguet and Joe Jackson and Nikiforos Mimikos-Stamatopoulos and Panagiotis E. Souganidis},
      year={2023},
      eprint={2312.11373},
      archivePrefix={arXiv},
      primaryClass={math.OC},
      url={https://arxiv.org/abs/2312.11373}, 
}

@article{cardaliaguet2023algebraic,
  title={An algebraic convergence rate for the optimal control of {McKean--Vlasov} dynamics},
  author={Cardaliaguet, Pierre and Daudin, Samuel and Jackson, Joe and Souganidis, Panagiotis E.},
  journal={SIAM Journal on Control and Optimization},
  volume={61},
  number={6},
  pages={3341--3369},
  year={2023},
  publisher={SIAM},
  doi={10.1137/22M1486789}
}

\end{document}